\documentclass[psamsfonts,11pt]{amsart}
 
\usepackage[usenames,dvipsnames]{color}
\usepackage{palatino}
\usepackage{eulervm}
\usepackage{xifthen}
\usepackage[utf8]{inputenc}
\usepackage{amsmath}
\usepackage{mathtools}
\usepackage{amsthm}
\usepackage{amsfonts}
\usepackage{amssymb}
\usepackage[margin=1in]{geometry}
\usepackage[all]{xy}
\usepackage[english]{babel}
\usepackage{graphicx}
\usepackage{xparse,zref-savepos}
\usepackage{tikz}
\usepackage{tikz-cd}
\usepackage{hyperref}
\usepackage[shortlabels,inline]{enumitem}
\usetikzlibrary{matrix,arrows,decorations.pathmorphing,shapes.misc}
\usepackage{capt-of}
\usepackage{appendix}
\usepackage{bbm}
\numberwithin{equation}{section}
\usepackage{courier}
\usepackage[capitalize,noabbrev]{cleveref}
\usepackage{pdfsync}
\newtheorem{theorem}{Theorem}[section]
\newtheorem{lemma}[theorem]{Lemma}
\newtheorem{proposition}[theorem]{Proposition}

\newtheorem{corollary}[theorem]{Corollary}
\theoremstyle{definition}
\newtheorem{example}[theorem]{Example}
\newtheorem{examples}[theorem]{Examples}
\newtheorem{definition}[theorem]{Definition}
\newtheorem{remark}[theorem]{Remark}
\newtheorem{remarks}[theorem]{Remarks}

\newcommand{\eps}{\varepsilon}

\newcommand{\im}{\textup{im} }

\newcommand{\R}{\mathbb{R}}
\newcommand{\N}{\mathbb{N}}
\newcommand{\Z}{\mathbb{Z}}

\newcommand{\lip}{\mathbf{Lip}}

\newcommand{\lippairs}{\mathbf{Lip}_\mathbf{pairs}}

\newcommand{\pmeas}{\mathcal{M}^+}
\newcommand{\fpmeas}{\mathcal{M}_{\textup{fin}}^+}
\newcommand{\fmeas}{\mathcal{M}_\textup{fin}}

\newcommand{\lipnorm}[1]{\|#1\|_{\textup{Lip}}}
\newcommand{\elp}[2]{\left\|#1\right\|_{#2}}
\newcommand{\cmon}[1][]{
	\ifthenelse{\isempty{#1}}{\mathbf{CMon}}{\mathbf{CMon}(#1)}
}
\newcommand{\ab}[1][]{
	\ifthenelse{\isempty{#1}}{\mathbf{Ab}}{\mathbf{Ab}(#1)}
}
\newcommand{\cmonti}[1][]{
	\ifthenelse{\isempty{#1}}{\mathbf{CMon}^{\mathbf{ti}}}{\mathbf{CMon}(#1)^{\mathbf{ti}}}
}
\newcommand{\abti}[1][]{
	\ifthenelse{\isempty{#1}}{\mathbf{Ab}^{\mathbf{ti}}}{\mathbf{Ab}(#1)^{\mathbf{ti}}}
}

\newcommand\restr[2]{{
		\left.\kern-\nulldelimiterspace 
		#1 
		\vphantom{\big|} 
		\right|_{#2} 
}}

\newcommand{\cat}{\mathbf}
\newcommand{\isom}{\cong}
\newcommand{\incl}{\hookrightarrow}

\DeclareMathOperator{\Cost}{Cost}
\DeclareMathOperator{\op}{op}
\DeclareMathOperator{\Lip}{Lip}
\DeclareMathOperator{\Intervals}{Int}

\makeatletter
\def\moverlay{\mathpalette\mov@rlay}
\def\mov@rlay#1#2{\leavevmode\vtop{%
		\baselineskip\z@skip \lineskiplimit-\maxdimen
		\ialign{\hfil$\m@th#1##$\hfil\cr#2\crcr}}}
\newcommand{\charfusion}[3][\mathord]{
	#1{\ifx#1\mathop\vphantom{#2}\fi
		\mathpalette\mov@rlay{#2\cr#3}
	}
	\ifx#1\mathop\expandafter\displaylimits\fi}
\makeatother

\DeclarePairedDelimiter{\abs}{\lvert}{\rvert}
\newcommand{\Expect}{{\rm I\kern-.3em E}}
\newcommand{\norm}[1]{\left\lVert#1\right\rVert}
\newcommand{\mW}{\mathbf{W}}

\makeatletter
\@ifundefined{zsaveposx}{\let\zsaveposx\zsavepos}{}
\newcounter{hposcnt}
\renewcommand*{\thehposcnt}{hpos\number\value{hposcnt}}
\NewDocumentCommand{\lplabel}{o m}{%
	\stepcounter{hposcnt}%
	\zsaveposx{\thehposcnt l}%
	\zref@refused{\thehposcnt l}%
	\zref@refused{hpos0l}%
	\makebox[0pt][r]{\makebox[\dimexpr\zposx{\thehposcnt l}sp-\zposx{hpos0l}sp][l]{#2}}%
	\IfNoValueF{#1}
	{\def\@currentlabel{#2}\ltx@label{#1}}
}
\makeatother
\AtBeginDocument{\zsaveposx{hpos0l}}

\newlist{thmenum}{enumerate}{1}
\setlist[thmenum, 1]{label=(\arabic*), ref=\thetheorem ~(\arabic*)}

\begin{document}
	
	\title{Virtual Persistence Diagrams, Signed Measures, Wasserstein Distances, and Banach Spaces}	
	\author{Peter Bubenik and Alex Elchesen}
	\maketitle
	
	
	\begin{abstract}
	Persistence diagrams, an important summary in topological data analysis, consist of a set of ordered pairs, each with positive multiplicity. Persistence diagrams are obtained via M\"obius inversion and may be compared using a one-parameter family of metrics called Wasserstein distances. In certain cases, M\"obius inversion produces sets of ordered pairs which may have negative multiplicity. We call these virtual persistence diagrams. Divol and Lacombe recently showed that there is a Wasserstein distance for Radon measures on the half plane of ordered pairs that generalizes both the Wasserstein distance for persistence diagrams and the classical Wasserstein distance from optimal transport theory. Following this work, we define compatible Wasserstein distances for persistence diagrams and Radon measures on arbitrary metric spaces. We show that the 1-Wasserstein distance extends to virtual persistence diagrams and to signed measures. In addition, we characterize the Cauchy completion of persistence diagrams with respect to the Wasserstein distances. We also give a universal construction of a Banach space with a 1-Wasserstein norm. Persistence diagrams with the 1-Wasserstein distance isometrically embed into this Banach space.
\end{abstract}
	
 	\section{Introduction}
	In computational settings, one-parameter persistent homology returns a finite indexed set of ordered pairs~\cite{edelsbrunner2000topological}, called a \emph{persistence diagram}. 
	The collection of persistence diagrams has a one-parameter family of metrics for $p \in [1,\infty]$ called \emph{Wasserstein distances}~\cite{cohen2007stability,cohen2010lipschitz}.
	The resulting metric spaces have a Cauchy completion~\cite{mileyko2011probability,blumberg2014robust}.
	Persistence diagrams and their Wasserstein distances are central to large parts of topological data analysis~\cite{Munch:2017,Robinson:2017,Seversky_2016_CVPR_Workshops}.
	
	More recently, our understanding of persistence diagrams and Wasserstein distance has been extended in the following ways.
	Persistence diagrams may be derived from the more elementary rank function via M\"obius inversion~\cite{Patel:2018}.
	The Wasserstein distance above may be extended to Radon measures on $\R^2_< = \{(x,y) \in \R^2 \ | \ x < y\}$ in a way that is compatible with the classical Wasserstein distance for probability measures~\cite{divol2019understanding}.
	By interpreting the collection of persistence diagrams algebraically as a free commutative monoid, the Wasserstein distances may be obtained in a functorial way, which implies that they have corresponding universal properties~\cite{bubenik2019universality}.
	The same construction may be applied to intervals, obtaining Wasserstein distances for barcodes~\cite{Collins:2004}, or to invariants of multiparameter persistence~\cite{bubenik2019universality}.
	If one applies M\"obius inversion to the graded rank function one obtains a variant of a persistence diagram in which the multiplicity of the ordered pairs is allowed to be negative~\cite{betthauser2019graded}. Generalized persistence diagrams which take on negative values have also arise naturally in several other places \cite{kim2021generalized, botnan2021signed, asashiba2019approximation, mccleary2020edit}. These ``virtual" persistence diagrams can be viewed as signed measures. Several extensions of the Wasserstein distance to the setting of signed measures have been introduced and studied~\cite{mainini2011}.

  The main goal of our work is to give various larger formal settings for generalized persistence diagrams (arising in one-parameter, multi-parameter, and generalized persistence) and their Wasserstein distances, which will be useful for new algorithms and theory in computational topology and the development of analytic tools for topological data analysis.
  As an example of the former, some of our results were used to define a stable Wasserstein distance for graded persistence diagrams~\cite{betthauser2019graded}.
  Our secondary goal is to help connect topological data analysis with optimal transportation.

\subsection{Our contributions}  
	
	We develop a number of constructions that unify and extend the results discussed above in various ways.
	All of our constructions are \emph{universal}. That is, they are functorial constructions satisfying certain universal properties. They may be interpreted as the existence of certain adjoint functors. Inspired by \cite{divol2019understanding}, instead of restricting ourselves to the usual setting, $\R^2_<$, we work in the general setting of metric spaces, to facilitate interactions with optimal transport theory and to permit applications to metric spaces of invariants for multiparameter persistent homology.
	
Let $(X,d)$ be a metric space and $A \subset X$. We call $(X,d,A)$ a \emph{metric pair} (see \cref{sec:metric-spaces}).
Formal sums on $(X,d,A)$ are (generalized) persistence diagrams~\cite{bubenik2019universality}.
  The following examples arise from persistent homology: consider $(\R^2_{\leq},d,\Delta)$, where $\R^2_{\leq} = \{(x,y) \in \R^2 \ | \ x \leq y\}$, $d$ is some metric on $\R^2_{\leq}$, and $\Delta = \{(x,y) \in \R^2 \ | \ x=y\}$;
replacing real numbers with extended real numbers, we have $(\overline{\R}^2_{\leq},d,\overline{\Delta})$; and
$(\Intervals(\R),d,\{\emptyset\})$, where $\Intervals(\R)$ denotes the set of intervals in $\R$, $d$ is some metric on this set, and $\emptyset$ denotes the empty interval.
Persistence diagrams and barcodes are formal sums on these metric pairs (see \cref{sec:persistence-diagrams}).
More generally, we have $(\Intervals(P),d,\{\emptyset\})$, where $P$ is some poset.
In particular, multiparameter persistent homology has $P=\R^d$ with the coordinatewise partial order. In this case, we can take $d$ to be Hausdorff distance or the volume (i.e. Lebesgue measure) of the symmetric difference~\cite{bubenik2018wasserstein}.
Most computational approaches to multiparameter persistent homology reduce to one-dimensional `slices'.
Let $\mathcal{L}$ be a set of lines (or, more generally, curves) in $\R^d$ that are images of order preserving maps $(\R,\leq) \to (\R^d,\leq)$.
Let $d'$ be a metric on $\R^2_{\leq}$ and $\{c_{\ell}\}_{\ell \in \mathcal{L}}$ be a set of nonnegative scaling constants.
If $\mathcal{L}$ is finite,
a set of persistence diagrams indexed by $\mathcal{L}$ is a formal sum in the metric pair $(\R^2_{\leq}\times \mathcal{L},d,\Delta \times \mathcal{L})$, where
the metric $d$ is given by $d((a,\ell),(b',\ell')) = c_{\ell}d'(a,a')$ if $\ell = \ell'$ and $\infty$ otherwise.
If $\mathcal{L}$ is finite or infinite, then given two sets of persistence diagrams indexed by $\mathcal{L}$ we may compute the Wasserstein distance for each $\ell \in \mathcal{L}$ and then compute the $p$-norm of the resulting function on $\mathcal{L}$ by summing or integrating with respect to an appropriate measure on $\mathcal{L}$.

We now list our main results.

	
	
	
	
	\subsubsection{Virtual persistence diagrams on metric pairs}
	Motivated by the growing number of settings in which signed persistence diagrams and signed barcodes arise \cite{betthauser2019graded, kim2021generalized, botnan2021signed, asashiba2019approximation, mccleary2020edit} and with the goal of systematically extending the Wasserstein distances to this setting, we develop a general framework for studying distances in the signed setting.
	Given a set $X$, let $K(X)$ denote the free abelian group on $X$.
	Let $(X,d,A)$ be a metric pair.
	Let $K(X,A)$ denote the quotient group $K(X)/K(A)$, which is isomorphic to $K(X \setminus A)$.
	We call the elements of $K(X,A)$ \emph{virtual persistence diagrams on $(X,A)$}.
	We prove there is a universal construction of the abelian group $(K(X,A),+)$ together with a metric $W_1$ given by
	\[
	W_1(\alpha^+ - \alpha^-, \beta^+ - \beta^-) = W_1(\alpha^+ + \beta^-, \beta^+ + \alpha^-)
      \]
      (see \cref{def:vpd}).
	We show that this metric is $1$-subadditive and \emph{translation invariant}.
	That is, $W_1(\alpha + \gamma, \beta + \gamma) = W_1(\alpha,\beta)$ (see \cref{cor:vpd-w1}). We also consider the corresponding constructions for the \emph{$p$-Wasserstein distances}, which are $\ell^p$ versions of the $1$-Wasserstein distance and which are widely studied in optimal transport theory.
	We give the analogous constructions for $W_p$, $p\in (1,\infty]$, but they require $d$ to be a \emph{$p$-metric}.
	That is, $d(x,y) \leq \norm{(d(x,z),d(z,y))}_p$ for all $x,y,z \in X$ (see \cref{thm:virtual_persistence_diagrams}).
	 
	
	We note that the form of this metric was already introduced by Mainini in the context of studying extensions of the classical Wasserstein distances to the setting of signed Radon measures \cite{mainini2011}. Our contribution here is that we derive this metric as the universal extension of a translation invariant metric compatible with the monoid structure of the space of persistence diagrams. More generally, we consider \emph{commutative metric monoids} $(M,d,+,0)$ and state conditions to ensure that the metric $\rho$ on the Grothendieck group $G(M)$ given by $\rho(m^+ - m^-, n^+-n^-) := d(m^+ + n^-, n^+ + m^-)$ is the canonical extension to the signed setting.

        \subsubsection{Measures on metric pairs}
	
	Given a metric space $(X,d)$, let $\pmeas(X)$ denote the commutative monoid of all Radon measures on $X$.
	Let $(X,d,A)$ be a metric pair, with $A\subset X$ a Borel subset, and let $p \in [1,\infty]$.
	We define $\pmeas(X,A)$ to be the quotient monoid $\pmeas(X)/\pmeas(A)$, which is isomorphic to $\pmeas(X \setminus A)$.
	We call the elements of $\pmeas(X,A)$ \emph{Radon measures on $(X,A)$} (\cref{def:relative_measures}). There is a $p$-subadditive metric $W_p$ on $\pmeas(X,A)$ which we call the \emph{Wasserstein distance}. Let $\pmeas_p(X,A)$ denote the submonoid of measures that are \emph{$p$-finite}. That is, $\int_{X \setminus A}d(x,A)^p\, d\mu < \infty$. Then $W_p$ restricts to a $p$-subadditive metric on $(\pmeas_p(X,A),+)$ (\cref{def:wasserstein-radon-pairs}). We show that this metric agrees with that introduced in \cite{divol2019understanding} for measures on $\R^2_<$ (\cref{cor:agrees-dl}). By taking $A = \emptyset$, we recover the classical Wasserstein distances between measures of equal mass. If we consider persistence diagrams to be discrete measures on $(X,A)$, then we obtain an inclusion $(D(X,A),+)\hookrightarrow (\pmeas(X,A),+)$ and corresponding isometric embedding $(D(X,A),W_p)\hookrightarrow (\pmeas(X,A),W_p)$ (\cref{prop:grothendieck_of_quotient}).
	
\subsubsection{Signed measures}
	
Given a metric space $(X,d)$, let $\fpmeas(X)$ and $\fmeas(X)$ denote the commutative monoid of finite Radon measures on $X$ and the abelian group of finite signed Radon measures on $X$, respectively.
We show that the metric $W_1$ extends from $\fpmeas(X)$ to its Grothendieck group $\fmeas(X)$ (\cref{prop:groth_group_of_measures,def:signed_wasserstein}).
We use the transshipment formulation of the $1$-Wasserstein distance (\cref{def:transshipment}) to show that $W_1$ on $\fmeas(X)$ is a solution to the signed transportation problem (\cref{def:signed-transportation,thm:signed-transportation}).

	Let $(X,d,A)$ be a metric pair, with $A\subset X$ a Borel subset.
	Let $\fmeas(X,A)$ denote the quotient group $\fmeas(X)/\fmeas(A)$, which is isomorphic to $\fmeas(X \setminus A)$.
	We call the elements of $\fmeas(X,A)$ \emph{finite signed Radon measures on $(X,A)$} (\cref{def:relative_measures}).
We show that $\fmeas(X,A)$ is the Grothendieck group of $\fpmeas(X,A)$ (\cref{prop:grothendieck_of_quotient}).
We define a commutative monoid $\pmeas_p(X,A)$ of $p$-finite Radon measures on $(X,A)$ and define a $p$-Wasserstein distance $W_p$ on $\pmeas_p(X,A)$ (\cref{def:wasserstein-radon-pairs,rem:wasserstein-radon-pairs}).
  We show that for classical persistence diagrams this definition agrees with the one of Divol and Lacombe~\cite{divol2019understanding} (\cref{cor:agrees-dl}).
As a result the $1$-Wasserstein distance extends to $M_1(\R^2_{\leq},\Delta)$ (\cref{cor:extending_classical_persistence_wasserstein}).
	
	\subsubsection{Cauchy completion of persistence diagrams on metric pairs}
	
	Given a set $X$,
	let $\overline{D}(X)$ denote the commutative monoid of all formal countable sums in $X$.
	Let $(X,d,A)$ be a metric pair and $p \in [1,\infty]$.
	Let $\overline{D}(X,A)$ denote the quotient $\overline{D}(X)/\overline{D}(A)$, which is isomorphic to $\overline{D}(X \setminus A)$.
	We call elements of $\overline{D}(X,A)$ \emph{countable persistence diagrams in $(X,A)$} (\cref{def:overline-d}).
	Let $\overline{D}_p(X,A)$ denote the submonoid of $\overline{D}(X,A)$ consisting of those countable persistence diagrams that, after removing at most finitely many summands, have finite $W_p$-distance to the zero persistence diagram, i.e., the \emph{$p$-finite} diagrams (\cref{def:diagrams_with_finite_distance_to_A}).
	We show that there is a universal construction of the commutative monoid $(\overline{D}_p(X,A),+)$ together with the metric $W_p$ such that $(\overline{D}_p(X,A),W_p)$ is a complete metric space 
        and $W_p$ is $p$-subadditive (\cref{thm:cauchy,thm:univ_prop_for_cauchy_completion_of_diagrams}).

\subsubsection{Free Banach spaces on metric pairs}

  Given a set $X$ there is a universal (real) vector space, $V(X)$, consisting of formal linear sums on $X$.
  Let $(X,d,A)$ be a metric pair. Let $V(X,A)$ denote the quotient $V(X)/V(A)$. We construct a universal normed vector space $(V(X,A),\norm{\ }_{W_1})$ (\cref{thm:universal-nvs}). Taking the Cauchy completion of this vector space, we obtain the free Banach space on the pair $(X,d,A))$ (\cref{thm:universal-banach}).

\subsubsection{Isometric embeddings}

Given a pointed metric space $(X,d,x_0)$ we have constructed a sequence of isometric embeddings (\cref{thm:isometric-embedding})
\begin{equation*} 
    (X,d) \incl (D(X,x_0),W_1) \incl (K(X,x_0),W_1) \incl (V(X,x_0),W_1) \incl (\hat{V}(X,x_0),W_1)
\end{equation*}
into metric free commutative monoid, a metric free abelian group, a normed vector space, and a Banach space, where each of the latter are canonically constructed.
Given a metric pair $(X,d,A)$ we have constructed a sequence of isometric embeddings (\cref{cor:isometric-embedding})
\begin{equation*} 
    (D(X,A),W_1) \incl (K(X,A),W_1) \incl (V(X,A),W_1) \incl (\hat{V}(X,A),W_1).
  \end{equation*}

\subsection{Remark on our constructions}
Our constructions are universal. That is, they arise as adjoints to certain forgetful functors. Benefits of having such constructions include the following.
\begin{enumerate*}
  \item We may be confident that our constructions are, in some sense, the right ones, rather than being ad-hoc.
  \item Our constructions have corresponding universal properties.
  \item Our constructions are functorial. Not only do we have various spaces of persistence diagrams but given any Lipschitz map between metric pairs we have a corresponding morphism between the resulting spaces of persistence diagrams. This may useful for metric pairs arising from multi-parameter and generalized persistence.
\end{enumerate*}

	\subsection{Related work}
	The Wasserstein distance has been thoroughly studied 
        in the context of measures
~\cite{rachev1998mass,villani2003topics}.
Persistence diagrams on $(\R^2_\leq,\Delta)$ and their Wasserstein distances were introduced in \cite{cohen2007stability,cohen2010lipschitz}.
  Persistence diagrams on $(\Intervals(\R),\emptyset)$ are called barcodes and were introduced in \cite{Collins:2004}.
Divol and Lacombe \cite{divol2019understanding}, connected the Wasserstein distances for measures and persistence diagrams on $(\R^2_{\leq},\Delta)$.
Turner and Skraba have 
given stability results for persistence diagrams on $(\R^2_\leq,\Delta)$ and their Wasserstein distances~\cite{Skraba:2020}.
%

The prequel to this paper~\cite{bubenik2019universality} introduces persistence diagrams on metric pairs and their Wasserstein distances. It shows that these arise from an adjoint functor on metric pairs and have corresponding universal properties. Furthermore, it is shown that the canonical inclusion of a metric pair into its space of persistence diagrams with the $p$-Wasserstein distance is $1$-Lipschitz and is an isometric embedding if $p=1$. In addition, it is shown that the $1$-Wasserstein distance satisfies Kantorovich-Rubenstein duality.

The free Banach space has been constructed and studied independently several times \cite{arens1956embedding, flood1984free, pestov1986free, godefroy2003lipschitz-free,weaver2018lipschitz}. Giusti and Lee~\cite{giusti2021signatures} have independently arrived at the connection between the Wasserstein distances and free Banach spaces (see Section~\ref{sec:universal-banach}).
The topological and metric properties of spaces of persistence diagrams have been studied by Che et al \cite{Che:2021} and Bubenik and Hartsock \cite{irynapeter}.

 	\subsection{Organization of the paper}

          In Section \ref{sec:background},we review background used throughout 
          and establish notation. In Section 3, we introduce a general theory of the Grothendieck group completion of monoids equipped with metric structures compatible with the monoid operation. We give sufficient conditions to guarantee that the metric extends from the monoid to its Grothendieck group completion in a compatible way. In Section \ref{sec:pers_wass},
we give a necessary and sufficient condition for the $p$-Wasserstein distance to be translation invariant, which may be of independent interest.          
We then introduce virtual persistence diagrams on metric pairs and define corresponding Wasserstein distances. In Section \ref{sec:measures}, we extend our constructions to Radon measures defined on metric pairs. In Section \ref{sec:cauchy_completion_of_pers_diags}, we extend our constructions to countably infinite persistence diagrams defined on metric pairs, and prove a universal property.
         In Section~\ref{sec:universal-banach}, we study universal constructions of vector spaces and Banach spaces.

\section{Background and notation}\label{sec:background}

In this section we summarize background material and corresponding notation that will be used throughout.

\subsection{Metric spaces and Lipschitz maps} \label{sec:metric-spaces}
\cite{bbi:book}.
	%
%
  A \emph{metric space} $(X,d)$ consists of a set $X$ and a function $d:X \times X \to [0,\infty]$ such that $d(x,x) = 0$ for all $x\in X$ (point equality), $d(x,y) = d(y,x)$ for all $x,y\in X$ (symmetry), and $d(x,y) \leq d(x,z) + d(z,y)$ for all $x,y,z\in X$ (triangle inequality).
Note that it is usually also assumed that a $d$ satisfies $d(x,x') = 0 \implies x = x'$ (separation) and $\im(d)\subset [0,\infty)$ (finiteness), but we will only make these assumptions in Section~\ref{sec:universal-banach}.
%
%
 With our definition, the $q$-norm for $q \in [1,\infty]$ induces not only a metric on the plane $\R^2$ but also on the extended plane $\overline{\R}^2$, where one doesn't have finiteness.
  Also, $\Intervals(\R)$ with either the Hausdorff distance or the length of the symmetric difference is a metric space which does not satisfy separation or finiteness.

	
	A function $f:X\to Y$ between metric spaces $(X,d_X)$ and $(Y,d_Y)$ is \textit{Lipschitz} if there exists a constant $C\geq 0$ such that
	$d_Y(f(x),f(x'))\leq C d_X(x,x')$
	for all $x,x'\in X$. In this case, $f$ is said to be \textit{$C$-Lipschitz} and $C$ is called a \textit{Lipschitz constant for $f$}. For $f:X\to Y$ Lipschitz, the \textit{Lipschitz norm} of $f$, denoted $\lipnorm{f}$, is defined
        to be the infimum of all such $C$.
	This infimum is in fact a minimum so that $d_Y(f(x),f(x'))\leq \lipnorm{f} d_X(x,x')$ for all $x,x'\in X$. If $f:X\to Y$, $g:Y\to Z$ are Lipschitz then $\lipnorm{g\circ f} \leq \lipnorm{f}\lipnorm{g}$.
	The collection of metric spaces together with Lipschitz maps forms a category denoted $\lip$.
%
		Metrics $d,\rho$ defined on the same set $X$ are said to be \textit{equivalent} if there are constants $C,K> 0$ such that
		$C\rho(x,x')\leq d(x,x')\leq K\rho(x,x')$
		for all $x,x'\in X$.

                A \emph{pair} is a tuple $(X,A)$ where $X$ is a set and $A\subset X$. A \textit{metric pair} is a tuple $(X,d,A)$ where $(X,d)$ is a metric space and $A \subset X$.
                When $A = \{x_0\}$ is a point then a pair $(X,\{x_0\})$ is also called a \textit{pointed set} and is denoted by $(X,x_0)$. Similarly, the metric pair $(X,d,\{x_0\})$ is called a \textit{pointed metric space} and is denoted $(X,d,x_0)$. 
                The distinguished point $x_0$ is called the \textit{basepoint}. A \textit{map of pairs} $f:(X,A)\to (Y,B)$ is a map $f:X\to Y$ such that $f(A) \subset B$. If $A  = \{x_0\}$ and $B = \{y_0\}$ then
                $f$ is a \textit{basepoint-preserving} (or \textit{pointed}) map.
          Metric pairs together with Lipschitz maps between metric pairs form a category which we denote by $\lippairs$. By $\lip_*$ we denote the full subcategory whose objects are pointed metric spaces.
%
		Let $X$ be a set, $(Y,d)$ a metric space, and $f:X\to Y$ any function. The \emph{pullback of $d$ through $f$} is a metric on $X$ defined by
		$f^*d(x,x') = d(f(x),f(x'))$
		for all $x,x'\in X$.
%
		A  Lipschitz map of pairs $f:(X,d_X,A)\to (Y,d_Y,B)$ is an isomorphism in $\lippairs$ if and only if $f$ is bijective, $f(A) = B$, and $d_X$ and $f^*d_Y$ are equivalent metrics on $X$.
In particular,		if $d$, $\rho$ are equivalent metrics on the metric pair $(X,A)$ then $(X,d,A)$ and $(X,\rho,A)$ are isomorphic objects in $\lippairs$.
%
	We denote the one-point metric space, viewed as a pointed metric space, by $*$. It is straightforward to check that $*$ is the terminal object in $\lippairs$, i.e., for every metric pair $(X,d_X,A)$ there exists precisely one Lipschitz map of pairs $f:X\to *$. $*$ is also the terminal object of the subcategory $\lip_*$.

	Let $(X,d,A)$ be a metric pair and let $X/A= (X\setminus A)\cup \{A\}$ denote the quotient set obtained by collapsing $A$ to a point. For $p\in [1,\infty]$, we define a metric $\overline{d}_p$ on $X/A$ by
	\[ \overline{d}_p(x,y)= \min\left(d(x,y),\elp{(d(x,A),d(y,A))}{p}\right).\]
	The proof that each $\overline{d}_p$ is indeed a metric can be found in \cite{bubenik2019universality}. Moreover, $\overline{d}_p$ metrizes the quotient topology on $X/A$ and satisfies the following universal property. Fix $p\in[1,\infty]$, let $(Y,d_Y,y_0)\in \lip_*$, and let $\phi:(X,d_X,A)\to (Y,d_Y,y_0)$ be a morphism in $\lippairs$. Then there is a unique basepoint preserving Lipschitz map $\tilde{\phi}:(X/A,\overline{d}_p,A)\to (Y,d_Y,y_0)$ satisfying $\tilde{\phi}\circ \pi = \phi$, where $\pi:X\to X/A$ denotes the quotient map.
        
        \begin{definition}[\cite{bubenik2019universality}]
          \label{def:product_metric}
  Let $X = (X,d_X)$ and $Y = (Y,d_Y)$ be metric spaces.
  For each $p\in [1,\infty]$ let $d\times_p d :X\times Y\to [0,\infty]$ be the function defined by
		\[ d\times_p d((x,y),(x',y'))= \elp{\left(d_X(x,x'),d_Y(y,y')\right)}{p},\]
		for all $x,x'\in X$, $y,y'\in Y$. We refer to $d\times_p d$ as the \emph{$p$-product metric on $X\times Y$}.
	\end{definition}
	
	We will also denote $d\times_1 d$ more simply by $d+d$. It follows from the Minkowski inequality that $d\times_p d$ is a metric on $X\times Y$ for each $p\in [1,\infty]$, with $d\times_p d$ and $d\times_q d$ being equivalent for any $p,q$. 
	Since $d\times_\infty d$ metrizes the product topology on $X\times Y$, it follows from metric equivalence that $d\times_p d$ also metrizes the product topology on $X\times Y$ for all $1\leq p\leq \infty$ \cite{bubenik2019universality}.
	Note that a metric $d:X\times X\to \R$ is Lipschitz
        with respect to the $p$-product metric.
	
		For any choice of $p\in [1,\infty]$, the canonical projections $\pi_X:X\times Y\to X$ and $\pi_Y:X\times Y \to Y$ are Lipschitz whenever $X\times Y$ is equipped with $d\times_p d$. When considering the product of a space $X$ with itself, we will denote the projections more simply by $\pi_1$ and $\pi_2$.
	\begin{proposition}\label{prop:prod_in_lip}
		Let $X = (X,d_X,A)$, $Y = (Y,d_Y,B)$ be metric pairs viewed as objects in $\lippairs$. For any choice of $p\in [1,\infty]$, the categorical product of $X$ and $Y$ in $\lippairs$ is given (up to isomorphism in $\lippairs)$ by $(X\times Y, d\times_p d, A\times B)$.

	\end{proposition}
	\subsection{Monoids}(\cite[I.1]{hungerford_algebra})
	A \textit{monoid} is a tuple $(M,+,0)$, where $M$ is a set, $+:M\times M\to M$ is an associative binary operation, and $0\in M$ is an identity element satisfying $0+m = m+0 = m$ for all $m\in M$. A monoid is \emph{commutative} if $m+n = n+m$ for all $m,n\in M$ and is \emph{cancellative} if $m + p = n + p \implies m = n$ for all $m,n,p\in M$. Submonoids are defined in the obvious way. \emph{Groups} are monoids in which every element has a two-sided inverse and are automatically cancellative. A \textit{monoid homomorphism} between monoids $M = (M,+_M,0_M)$ and $N = (N,+_N,0_N)$ is a function $f:M\to N$ such that $f(m+_M n) = f(m) +_N f(n)$ for all $m,n\in N$ and $f(0_M) = 0_N$. We remark that, unlike for group homomorphisms, the requirement that $f(0_M) = 0_N$ does not follow automatically from the first condition.

	The \emph{free commutative monoid} on a nonempty set $X$, denoted $D(X)$, is defined by
	\[D(X) = \{f:X\to \N \ | \ f(x) = 0 \textup{ for all but finitely many $x\in X$}\},\]
	with the monoid operation being addition of functions. We also set $D(\emptyset) = 0$, the zero monoid. The \emph{indicator for $x\in X$} is defined by $1_x(x') = 1$ if $x' = x$ and is $0$ otherwise. Note that every element of $D(X)$ can be written uniquely as a finite sum of such indicators. By identifying the indicator $1_x$ with $x$, we identify elements of $D(X)$ with (finite) formal sums of elements of $X$, with addition of functions corresponding to addition of formal sums. Since $\N$ is cancellative, so is $D(X)$. If $A\subset X$ then $D(A)$ can be identified as a submonoid of $D(X)$ in the obvious way. For a given set $X$, there is a canonical map $i = i_X:X\to D(X)$ given by $i(x) = x$ (the right-hand-side being a formal sum with just one summand). The free commutative monoid together with this canonical map satisfy the following universal property. 

\begin{proposition}\label{prop:free_comm_monoid_universality}
  For any commutative monoid $M$ and map $\phi:X\to M$, there exists a
  unique monoid homomorphism $\tilde{\phi}:D(X)\to M$ such that
  $\phi = \tilde{\phi}\circ i$.
\end{proposition}

If $h:X\to Y$ is any function then we define $Dh:D(X)\to D(Y)$ by
$Dh(x_1 + \dots + x_n) = h(x_1) + \dots + h(x_n)$.
Equivalently, viewing elements of $D(X)$ as functions, for $f:X\to \N$ with finite support, we have $Dh(f) = h_*(f)$ where for $y \in Y$, $h_*(f)(y) = \sum_{x \in h^{-1}(y)}f(x)$.
Then $Dh$ is a monoid homomorphism. Moreover, the assignments $X\mapsto D(X)$, $h\mapsto Dh$ specify a functor $D:\mathbf{Set}\to \cmon$, where $\mathbf{Set}$ denotes the category of sets and functions and $\mathbf{CMon}$ denotes the category of commutative monoids and monoid homomorphisms. Letting $U:\cmon\to \mathbf{Set}$ denote the forgetful functor, Proposition \ref{prop:free_comm_monoid_universality} is equivalent to the statement that $D$ is left adjoint to $U$.
	
	An equivalence relation $\sim$ on a monoid $M$ is called a \emph{congruence} if $a\sim b$ and $c\sim d$ implies $a+c \sim b + d$. If $\sim$ is a congruence then there is a well-defined monoid structure on the set of equivalence classes $M/\sim$ defined by $[a] + [b] = [a+b]$. 
%
		Let $M$ be a commutative monoid and let $N\subseteq M$ be any submonoid. Define a relation $\sim$ on $M$ by 
		$a\sim b$ iff $\exists x,y\in N$ such that $a+x = b+y$.
		It is easily verified that $\sim$ is a congruence. We denote the commutative monoid $M/\sim$ by $M/N$ and refer to it as the \textit{quotient of $M$ by $N$} and the equivalence class of $a\in M$ under the congruence $\sim$ is denoted $a+N$. We will denote the congruence of elements $a$ and $b$ under the above congruence relation by $a = b \pmod{N}$ and say that \emph{$a$ equals $b$ mod $N$}. Note that if $M$ is cancellative then so is $M/N$.
	
	Of particular importance for us is the quotient $D(X)/D(A)$, where $A\subset X$. In this case, $f=g \pmod{D(A)}$ if and only if $f|_{X \setminus A} = g|_{X \setminus A}$.
        Equivalently,
	$x_1 + \dots + x_m = x_1' + \dots + x_n' \pmod{D(A)}$ iff $x_1 + \dots + x_m + a_1 + \dots + a_s = x_1' + \dots + x_n' + a_1' + \dots + a_t'$,
	for some $a_1,\dots,a_s,a_1',\dots a_t'\in A$. Hence, $D(X)/D(A)\cong D(X\setminus A)$.
%
	If $f:(X,A)\to (Y,B)$ is a map of pairs then the \emph{induced map} $Df:D(X,A)\to D(Y,B)$ is defined by $x_1 + \dots + x_n + D(A)\mapsto f(x_1) + \dots + f(x_n) + D(B)$.

\subsection{Persistence diagrams}{\cite{bubenik2019universality}}.
\label{sec:persistence-diagrams}        
%
%
Let $(X,A)$ be a pair. A \emph{persistence diagram} on $(X,A)$ is an element of the commutative monoid $D(X,A) = D(X)/D(A)$.
%
          The monoid $D(\R^2,\R^2_\geq) \isom D(\R^2_{\leq},\Delta)$ is the monoid of classical persistence diagrams with finitely many points. 
The monoid $D(\Intervals(R),\emptyset)$ is the monoid of barcodes.

We remark that other terminology may be preferable for some readers. Persistence diagrams on a metric pair could instead be referred to as formal sums on a metric pair or discrete measures on a metric pair.

\subsection{Wasserstein distance}{\cite{bubenik2019universality}}.
%
%
%
%
		Let $(X,d,A)$ be a metric pair and let $\alpha,\beta\in D(X,A)$. Let $\pi_1,\pi_2:X\times X\to X$ denote the canonical projects. We call $\sigma = (x_1,x_1') + \dots + (x_r,x_r')\in D(X\times X)$ a \emph{matching between $\alpha$ and $\beta$} if $(\pi_1)_*\sigma = \alpha \pmod{D(A)}$ and $(\pi_2)_*\sigma = \beta \pmod{D(A)}$. Let $p\in [1,\infty]$. The \emph{$p$-cost of $\sigma$ (with respect to $d$)} is defined by
		\[\Cost_p[d](\sigma) = \elp{\big(d(x_i,x_i')\big)_{i = 1}^r}{p}.\]
%
                Let $(X,d,A)$ be a metric pair. For each $p\in[1,\infty]$ the \emph{$p$-Wasserstein distance} $W_p[d]:D(X,A)\times D(X,A)\to [0,\infty]$ is defined as follows. Given $\alpha,\beta\in D(X,A)$, we set
		\[ W_p[d](\alpha,\beta) = \inf_{\sigma}\,\textup{Cost}_p[d](\sigma),\]
		the infimum being taken over all matchings between $\alpha$ and $\beta$. 
		
				\begin{remark}	For $p \in [1,\infty)$, we could equivalently define the $p$-Wasserstein distance by defining $\textup{Cost}_p[d](\sigma) = {\textstyle\sum_{i = 1}^r d(x_i,x_i')^p}$ and then setting $W_p[d](\alpha,\beta) = \left(\inf_{\sigma}\textup{Cost}_p[d](\sigma)\right)^{1/p}$. This is the convention often used in the optimal transport literature \cite{villani2003topics}. We will stick to the convention introduced above for persistence diagrams as it includes the $p = \infty$ case as well, but use the optimal transport convention later when we work with measures (see Definitions \ref{def:transportation_problem}, \ref{def:wasserstein-radon-pairs}, and \ref{def:partial_wass_divol_lacombe}).
	\end{remark}
	
	
                Consider the metric pair $(\R^2_\leq,d,\Delta)$ where $d$ is the metric induced by the $\infty$-norm and the monoid $D(\R^2_\leq,\Delta)$ of classical persistence diagrams.
                Let $\alpha = (0,1) + (2,4)$ and let $\beta = (3,5)$. Then $\sigma = ((0,1),(1/2,1/2)) + ((2,4),(3,5))$ is matching between $\alpha$ and $\beta$ with $C_p[d](\sigma) = 3/2$. Thus $W_p[d](\alpha,\beta)\leq 3/2$, and it is not too hard to see that this matching is optimal.
	
                        Let $p\in [1,\infty]$. A metric $d$ on a monoid $(M,+)$ is called \textit{$p$-subadditive}~\cite{bubenik2018wasserstein} if 
		\[ d(a+b,a'+b')\leq \elp{\left(d(a,a'), d(b,b')\right)}{p},\]
		for all	$a,a',b,b'\in M$. 
%
%
	The metric $W_p[d]$ is $p$-subadditive for any metric pair $(X,d,A)$ and $p\in [1,\infty]$ \cite[Lemma 4.17]{bubenik2019universality}. 

\subsection{The Grothendieck group completion}{\cite[II.1]{Weibel:kbook}}. \label{sec:grothendieck-group}
	Monoids can be completed into groups in a universal way. Given a commutative monoid $M = (M,+)$, the \emph{Grothendieck group completion of $M$} (or just the \emph{Grothendieck group of $M$}), denoted $K(M)$, is defined as follows. Let $\sim$ denote the equivalence relation on $M\times M$ given by $(a,b)\sim (c,d)$ if and only if there exists $k\in M$ such that $a+d + k = b + c + k$. As a set, we define $K(M) = (M \times M)/\sim$, the set of equivalence classes of $M\times M$ under $\sim$. We denote the equivalence class of $(a,b)$ by $a-b$. The binary operation on $K(M)$ (also denoted $+$) is given by 
	\[(a-b) + (c-d) = (a+c) - (b+d).\]
	It is straightforward to check that this operation is well-defined and makes $K(M)$ into an abelian group with identity $0 = 0-0$ and with inverses $-(a-b) = b-a$ for all $a-b\in K(M)$. We caution that $a-b = c-d$ does not necessarily imply that $a+d = b+c$. However, if $M$ is cancellative then $a-b = c-d$ iff $a+d = b+c$.
        For example, the Grothendieck group of the monoid $\N$ of natural numbers is isomorphic to the group $\Z$ of integers.
	
	The canonical map $u:M\to K(M)$ is the monoid homomorphism given by $u(a) = a - 0.$ If we use the notation $a= a-0$ in $K(M)$ then the canonical map takes the form $u(a) = a$ for all $a\in M$. Again, we must caution that with this notation, $a = b$ in $K(M)$ if and only if there exists some $k\in M$ such that $a + k = b + k$ in $M$, but that when $M$ is cancellative we have $a = b$ in $K(M)$ iff $a = b$ in $M$. Since all of the monoids we encounter in this paper will be cancellative, this convention should cause no confusion.
	The Grothendieck group together with the canonical map $u$ satisfy the following universal property.
	
	\begin{proposition}\label{prop:grothendieck_universalitiy}
		For any abelian group $H$ and monoid homomorphism $v:M\to H$, there exists a unique group homomorphism $\phi:K(M)\to H$ such that $\phi\circ u = v$.
	\end{proposition}

	For a monoid homomorphism $f:M\to N$ between commutative monoids $M,N$, the function $Kf:K(M)\to K(N)$ given by $a-b\mapsto f(a) - f(b)$ is a well-defined group homomorphism, and this assignment makes $G$ into a functor $K:\cmon\to \ab$. Let $U:\ab\to\cmon$ denote the inclusion functor of abelian groups into commutative monoids. Then Proposition \ref{prop:grothendieck_universalitiy} is equivalent to the statement that $K$ is left adjoint to $U$.
	
		For cancellative monoids, the description of the Grothendieck group of a quotient monoid is particularly simple.
		
		\begin{proposition}\label{prop:groth_of_quotient}
			Let $M$ be a cancellative commutative monoid and let $N$ be a submonoid of $M$. Then $K(M/N)\cong K(M)/K(N)$.
			\begin{proof}
				First, note that since $M$ is cancellative, $K(N)$ can be identified naturally as a subgroup of $K(M)$. Now define $\phi:K(M)\to K(M/N)$ by $m-m'\mapsto (m+N)-(m'+N)$. Clearly $\phi$ is surjective, and since $M/N$ is cancellative, $m-m'\in\textup{ker}(\phi)\iff m+N = m'+N\iff m+n=m'+n'$ for some $n,n'\in N$. But now $m+n = m'+n'\iff m-m' = n-n'\in K(N)$ so that $\textup{ker}(\phi) = K(N)$, and the result now follows from the first isomorphism theorem.
			\end{proof}
                   \end{proposition}

                   \subsection{Universal constructions}
                   \label{sec:universal}
                   
                   {\cite[Sections 2.3, 2.4, 4.2, 4.6]{Riehl:2016}}.
Many of our main results are stated in terms of the existence of an object with a property called a universal property. Equivalently, they may be stated as the existence of a certain adjoint functor. For a discussion of universal properties and adjunctions including the example of a free commutative monoid, see \cite[Section 2.5]{bubenik2019universality}.
	
	\section{Grothendieck groups of Lipschitz monoids}\label{sec:grothendieck_lipschitz_group}
	Many objects carry both the structure of a metric and of a commutative monoid. Persistence diagrams, which motivate this work, are an example, as we will see. It is natural to ask whether or not the metric structure of a given monoid can be extended to its Grothendieck group in a natural way. The purpose of this section is to describe sufficient 
	conditions on the metric that guarantee this is possible.
	
\subsection{Lipschitz monoids and Lipschitz groups}
	
A \emph{commutative Lipschitz monoid} (CL monoid) is a commutative monoid internal to $\lip$. Explicitly, this is a metric space $(M,d)$ equipped with a binary operation $+$ and an identity element $0$ so that $(M,+,0)$ is a commutative monoid and $+:M\times M\to M$ is Lipschitz, where $M\times M$ is the categorical product in $\cat{Lip}$ (see \cref{prop:prod_in_lip}). Using the metric $d\times_p d$ on $M\times M$, 
the latter condition means that
	\[ d(x+y,x'+y')\leq \lipnorm{+}(d(x,x,') + d(y,y')),\]
	for all $x,x',y,y'\in M$. A morphism between CL monoids is a Lipschitz monoid homomorphism. CL monoids together with these morphisms form the category $\cmon[\lip]$.
	    
	An \emph{abelian Lipschitz group} (AL group) is an abelian group internal to $\lip$. This requires, in addition to the requirements for a CL monoid, that the map $-:M\to M$ taking a group element to its inverse be Lipschitz. AL groups together with Lipschitz group homomorphisms form the category $\ab[\lip]$.
	
	We will also consider the category $\cmon[\lippairs]$ of metric pairs with compatible monoid structures, the morphisms being Lipschitz maps of pairs which are also monoid homomorphisms. Similarly, we have the category $\cmon[\lip_*]$ of pointed metric spaces with compatible monoid structures, whose morphisms are basepoint-preserving Lipschitz monoid homomorphisms.
	\subsection{Translation invariant metrics on monoids}
	A metric $d$ on a commutative monoid $(M,+)$ is \emph{translation invariant} if $d(m+p, n+p) = d(m,n)$ for all $m,n,p\in M$. 
	
	
	\begin{proposition}\label{prop:ti_implies_CL_or_AL}
		Let $(M,+)$ be a commutative monoid equipped with a translation invariant metric $d$. Then
		\begin{enumerate}
			\item $(M,d,+)$ is a CL monoid with $\lipnorm{+}= 1$, unless $d = 0$ in which case $\lipnorm{+} = 0$; and
			\item If $M$ is a group then $(M,d,+)$ is an AL group.
		\end{enumerate}
		\begin{proof}
			
			(1) If $d = 0$ then the conclusion is immediate, so suppose that $d \neq 0$. By the triangle inequality and translation invariance for $d$ we have
			\[d(x+y,x'+y')\leq d(x+y,x'+y) + d(x'+y,x'+y') = d(x,x') + d(y,y').\]
			Thus the monoid operation is Lipschitz with $\lipnorm{+}\leq 1$. On the other hand, if $a,b\in M$ are such that $d(a,b)\neq 0$ and $c\in M$ is any other element then we have $d(a+c,b+c) = d(a,b) = d(a,b) + d(c,c)$, from which it follows that $\lipnorm{+}\geq 1$ as well.
			
			(2) By part (1), $+$ is Lipschitz. By translation invariance of $d$ we have $d(-a,-b) = d(a,b)$. Thus the inversion operation is an isometry, and hence Lipschitz.
		\end{proof}	
	\end{proposition}

	We will also consider the following weakening of the cancellative property.
	
	\begin{definition}
		A commutative monoid $M$ equipped with a metric $d$ is said to be \emph{weakly cancellative} if $a+c = b+c \implies d(a,b) = 0$ for all $a,b,c\in M$.
	\end{definition}

	\begin{lemma}\label{lem:ti_implies_weakly_cancellative}
		If $M$ is a commutative monoid equipped with a translation invariant metric $d$ then $M$ is weakly cancellative. Moreover, if $d$ satisfies the separation axiom then $M$ is cancellative.
		\begin{proof}
			If $a,b,c\in M$ are such that $a+c = b+c$ then $0 = d(a+c,b+c) = d(a,b)$ so that $M$ is weakly cancellative. If $d$ satisfies the separation axiom then $d(a,b) = 0 \implies a = b$ so that $M$ is cancellative.
		\end{proof}
	\end{lemma}
	
	\begin{corollary}\label{cor:equality_in_Groth_implies_weak_equality}
		If $M$ is a commutative monoid equipped with a translation invariant metric $d$ and $a-b = a'-b'$ in $K(M)$ then $d(a+b', a'+b) = 0$. In particular, if $ a =b$ in $K(M)$ then $d(a,b) = 0$.
		\begin{proof}
			Since $a-b = a'-b'$ in $K(M)$, there exists some $c\in M$ such that $a+b' + c = a' + b + c$ in $M$. By Lemma \ref{lem:ti_implies_weakly_cancellative}, $M$ is weakly cancellative and hence $d(a+b',a'+b) = 0$. The second statement follows by letting $a' = b' = 0$.
		\end{proof}
	\end{corollary}

	If $M$ is only weakly cancellative then the canonical map $u:M\to K(M)$ need not be injective. However, by also weakening the definition of injectivity we obtain something analogous.
	\begin{definition}\label{def:weak_injectivity}
		Let $(X,d)$ be a metric space and $Y$ any set. A map $f:X\to Y$ is \emph{weakly injective} if $f(x) = f(x')\implies d(x,x') = 0$.
	\end{definition}
	
	\begin{proposition}\label{prop:weak_canc_iff_weak_inj}
		Let $M$ be a commutative monoid equipped with a metric $d$ and let $K(M)$ be the Grothendieck group of $M$. Then $M$ is weakly cancellative if and only if the canonical map $u:M\to K(M)$ is weakly injective.
		\begin{proof}
			Let $a,b,c\in M$. If $M$ is weakly cancellative and $u(a) = u(b)$ then $a+k = b+k$ for some $k\in M$. Hence $d(a,b) = 0$ by weak cancellation so that $u$ is weakly injective. Conversely, if $u$ is weakly injective and $a+c = b+c$ in $M$ then $u(a) + u(c) = u(a + c) =  u(b+c)= u(b) + u(c)$ so that $u(a) = u(b)$. Hence $d(a,b) = 0$ by weak injectivity so that $M$ is weakly cancellative.
		\end{proof}
	\end{proposition}
	
	We will denote the full subcategories of $\cmon[\lip]$ and $\ab[\lip]$ whose objects are translation invariant commutative Lipschitz monoids (TICL monoids) and translation invariant abelian Lipschitz groups (TIAL groups) by $\cmonti[\lip]$ and $\abti[\lip]$, respectively.

	\subsection{The Grothendieck Lipschitz group}
	There is an obvious inclusion functor $U:\abti[\lip]\to \cmonti[\lip]$. Our goal is to show that a metric can be put on the Grothendieck group of a TICL monoid so as to define a functor $K:\cmonti[\lip]\to \abti[\lip]$ that is left adjoint to $U$. 
	
	\begin{definition}\label{def:Grothendieck_Lipschitz_group}
		Given a TICL monoid $(M,d,+)$, let $(K(M),+)$ be the corresponding Grothendieck group completion of $(M,+)$. We define a function $\rho:K(M)\times K(M)\to [0,\infty]$ by
		\begin{equation}\label{eq:grothendieck_metric}
			\rho(a-b,c-e) = d(a+e,c+b).		
		\end{equation}
	\end{definition}
	
	We will show that \eqref{eq:grothendieck_metric} is a metric which makes $K(M)$ into a TIAL group. 
	
	\begin{proposition}\label{prop:grothendick_metric_welldefined+ti}
		Let $K(M)$ be the Grothendieck group of a TICL monoid $(M,d,+)$. Then the function $\rho$ given by \eqref{eq:grothendieck_metric} defines a translation invariant metric on $K(M)$.
		\begin{proof}
			First, we verify that $\rho$ is well-defined. Suppose that $a-b = a'-b'$ and $c-e = c'-e'$ in $K(M)$. Then $d(a+b',a'+b) = 0 = d(c+e',c'+e)$ by Corollary \ref{cor:equality_in_Groth_implies_weak_equality}. Hence, by translation invariance and the triangle inequality for $d$,
$				\rho(a-b,c-e) = d(a+e,c+b) = d(a+e+b',c+b+b') 
				\leq d(a+e+b',a'+b+e) + d(a'+b+e,c+b+b') 
				= d(a+b',a'+b) + d(a'+e,c+b') = d(a'+e,c+b')$,
			and furthermore,
$				d(a'+e,c+b') = d(a'+e+e',c+b'+e') 
				\leq d(a'+e+e',c'+b'+e) + d(c'+b'+e,c+b'+e') 
				= d(a'+e',c'+b') = \rho(a'-b',c'-e')$,
			since $d(c'+b'+e,c+b'+e') = d(c'+e,c+e') = 0$. Thus $\rho(a-b,c-e) \leq \rho(a'-b',c'-e')$. The reverse inequality is obtained symmetrically, which shows that $\rho$ is well-defined.
			
			To see that the point equality holds for $\rho$, note that if $a-b = a'-b'$ then, since $d$ is translation invariant, $\rho(a-b,a'-b') = d(a+b',a'+b) = 0$ by Corollary \ref{cor:equality_in_Groth_implies_weak_equality}. Symmetry of $\rho$ follows immediately from the symmetry of $d$. To verify the triangle inequality for $\rho$, let $a-b,c-e,x-y\in K(M)$. Then by the triangle inequality for $d$,
$				\rho(a-b,x-y) + \rho(x-y,c-e) = d(a+y,x+b) + d(x+e,c+y) 
				= d(a+y+e,x+b+e) + d(x+b+e,c+b+y) 
				\geq d(a+y+e,c+b+y) = d(a+e,c+b) = \rho(a-b,c-e)$.
			Thus $\rho$ is a metric on $K(M)$. 
			
			Lastly, by translation invariance of $d$ we have
$				\rho((a-b) + (x-y),(c-e) + (x-y)) = \rho((a+x) - (b+y),(c+x) - (e+y)) 
				= d(a+x+e+y,c+x+b+y) = d(a+e,c+b) = \rho(a-b,c-e)$,
			which shows that $\rho$ is translation invariant as well.
		\end{proof}
	\end{proposition}
	
	\begin{corollary}\label{cor:grothendieck_is_TIAL}
		$(K(M),\rho,+)$ is a TIAL group, i.e., an object of $\abti[\lip]$.
		\begin{proof}
			By Proposition \ref{prop:grothendick_metric_welldefined+ti}, $\rho$ is translation invariant and so by Proposition \ref{prop:ti_implies_CL_or_AL}, $(K(M),\rho,+)$ is an AL group.
		\end{proof}
	\end{corollary}
	
	\begin{definition}
		When the canonical map $u:M\to K(M)$ is weakly injective, we say that a metric $\rho$ on $K(M)$ \emph{extends $d$} if $u^*\rho = d$.
	\end{definition}
	
	\begin{proposition}\label{prop:uniqueness_of_groth_metric}
		Let $(M,d,+)$ be a TICL monoid and let $K(M)$ be its Grothendieck group. The metric $\rho$ defined by \eqref{eq:grothendieck_metric} is the unique translation invariant metric extending $d$ to $K(M)$.
		\begin{proof}
			Since $d$ is translation invariant, $u$ is weakly injective by Lemma \ref{lem:ti_implies_weakly_cancellative} and Proposition \ref{prop:weak_canc_iff_weak_inj}. Furthermore, $\rho$ is translation invariant by Proposition \ref{prop:grothendick_metric_welldefined+ti} and $\rho(u(a),u(b)) = \rho(a-0,b-0) = d(a,b)$. Thus $\rho$ is a translation invariant metric on $K(M)$ extending $d$. Moreover, if $\eta$ is any such metric then 
$				\eta(a-b,c-d) = \eta((a-b )+ (b+d),(c-d) + (b+d)) = \eta(a+d,b+c) 
				= \eta(u(a+d),u(b+c)) = d(a+d,b+c) = \rho(a-b,c-d)$. 
		\end{proof}
	\end{proposition}


	\begin{definition}
		Given a TICL monoid $M = (M,d,+)$, we call the TIAL group $(K(M),\rho,+)$ the \emph{Grothendieck Lipschitz group} of $M$.
	\end{definition}
	
	\begin{proposition}\label{prop:canonical_map_lipschitz} Let $(M,d,+)$ be a TICL monoid.
		\begin{thmenum}
			\item \label{prop:canonical_map_lipschitz_part_1}	The canonical map $u:(M,d,+)\to (K(M),\rho,+)$ is Lipschitz with $\lipnorm{u} = 1$ if $d \neq 0$ and $\lipnorm{u} = 0$ otherwise.
			\item \label{prop:canonical_map_lipschitz_part_2} $\lipnorm{+_{K(M)}} = \lipnorm{+_M}$.
		\end{thmenum}
		\begin{proof}
			(1) Clearly if $d = 0$ then $\lipnorm{u} = 0$. Suppose that $d \neq 0$.  Since $\rho(u(a),u(b)) = d(a,b)$ for all $a,b\in M$, $\lipnorm{u}\leq 1$. On the other hand, given $a,b\in M$ with $d(a',b') \neq 0$, we have $\rho(u(a),u(b))/d(a,b) = 1$ and hence $\lipnorm{u}\geq 1$ as well.
			
			(2) Note that $\rho = 0 \iff d = 0$. It follows then from Proposition \ref{prop:ti_implies_CL_or_AL} that $\lipnorm{+_{K(M)}} = \lipnorm{+_M} = 0$ if $d = 0$ and $\lipnorm{+_{K(M)}} = \lipnorm{+_M} = 1$ otherwise.
		\end{proof}
	\end{proposition}
	
	The Grothendieck Lipschitz group satisfies the following universal property.
	
	\begin{theorem}\label{thm:univ_of_groth_lip}
		Let $M = (M,d,+)$ be a TICL monoid and let $K(M) = (K(M),\rho,+)$ denote its Grothendieck Lipschitz group. 
		
		\begin{thmenum}
			\item \label{thm:univ_of_groth_lip_part_1}For any TIAL group $H$ and any Lipschitz monoid homomorphism $\phi:M\to H$, there exists a unique Lipschitz group homomorphism $\tilde{\phi}:K(M)\to H$ such that $\tilde{\phi}\circ u = \phi$. 
			\item \label{thm:univ_of_groth_lip_part_2} $\lipnorm{\tilde{\phi}} = \lipnorm{\phi}$.
		\end{thmenum}
		\begin{proof}
			(1) By the universal property of the Grothendieck group, there is a unique group homomorphism $\tilde{\phi}:K(M)\to H$ such that $\tilde{\phi}\circ u = \phi$. It remains only to check that $\tilde{\phi}$ is Lipschitz. Let $d_H$ denote the metric on $H$. Then
			\begin{multline*}
				d_H(\tilde{\phi}(a-b),\tilde{\phi}(a'-b')) = d_H(\tilde{\phi}(a)-\tilde{\phi}(b),\tilde{\phi}(a')-\tilde{\phi}(b'))\\
				= d_H(\tilde{\phi}(a) + \tilde{\phi}(b'),\tilde{\phi}(a') + \tilde{\phi}(b)) = d_H(\tilde{\phi}(a+b'),\tilde{\phi}(a'+b))\\
				= d_H(\tilde{\phi}(u(a+b')),\tilde{\phi}(u(a'+b))) = d_H(\phi(a+b'),\phi(a'+b))\\
				\leq \lipnorm{\phi}d(a+b',a'+b) = \lipnorm{\phi}\rho(a-b,a'-b').
			\end{multline*}
			Thus $\tilde{\phi}$ is Lipschitz with $\lipnorm{\tilde{\phi}}\leq \lipnorm{\phi}$.
			
			(2) By Proposition \ref{prop:canonical_map_lipschitz_part_1}, $\lipnorm{u} = 1$. Hence $\lipnorm{\phi} = \lipnorm{\tilde{\phi}\circ u} \leq \lipnorm{\tilde{\phi}}\lipnorm{u} = \lipnorm{\tilde{\phi}}$. By the proof of (1), $\lipnorm{\tilde{\phi}}\leq \lipnorm{\phi}$ and the result follows.
		\end{proof}
	\end{theorem}
	
	Consider again the Grothendieck group completion functor $K:\cmon\to \ab$. Given TICL monoids $(M,d_M,+)$, $(N,d_N,+)$, consider their corresponding TIAL groups $(K(M),\rho_M,+)$, $(K(N),\rho_N,+)$, and let $f:M\to N$ be a Lipschitz monoid homomorphism. Then for $a-b,a'-b'\in K(M)$ we have
	\begin{align*}
		\rho_N(Kf(a-b),Kf(a'-b')) &= \rho_N(f(a)-f(b), f(a') - f(b'))\\
		& = d_N(f(a) + f(b'),f(a') + f(b))\\
		& = d_N(f(a+b'),f(a'+b))\\
		& \leq \lipnorm{f}d_M(a+b',a'+b) = \lipnorm{f}\rho_M(a-b,a'-b').
	\end{align*} 
	Thus $Kf$ is a Lipschitz group homomorphism with $\lipnorm{Kf}\leq \lipnorm{f}$. Hence we obtain a functor $K:\cmonti[\lip]\to\abti[\lip]$. Then Theorem \ref{thm:univ_of_groth_lip} implies the following.
	
\begin{corollary}
		$K$ is left adjoint to the inclusion functor $U:\mathbf{Ab}(\mathbf{Lip})^{\mathbf{ti}}\to \mathbf{CMon}(\mathbf{Lip})^{\mathbf{ti}}$.
	\end{corollary}
	
	\begin{example}\label{ex:naturals_grothendieck_lipschitz_group_completion}
		Continuing with the example of the natural numbers (Section~\ref{sec:grothendieck-group}), the absolute value norm on $\N$ induces a translation invariant metric $d_\N$ given by
		\[ d_\N(n,m) = |n-m| = \max(n,m)- \min(n,m).\]
		By Proposition \ref{prop:uniqueness_of_groth_metric}, $d_\N$ extends to a unique translation invariant metric $\rho_\Z$ on $\Z$ given by
		\[ \rho_\Z(a-b,c-d) = d_\N(a+d,b+c) = \max(a+d,b+c) - \min(a+d,b+c).\]
		This is of course the usual metric on $\Z$ induced by absolute value, more readily seen by using the absolute value notation: $\rho_\Z(a-b,c-d) = |a+d - (b+c)| = |(a-b) - (c-d)|$.
		
	\end{example}

\section{Virtual persistence diagrams}\label{sec:pers_wass}

In this section, we define virtual persistence diagrams on metric pairs and define Wasserstein distances between these objects. 

	\subsection{Translation invariance of the Wasserstein distances}
	Since we are interested in translation invariant metrics, we would like to determine for which $p$ and for what conditions on $d$ we are assured that $W_p[d]$ is translation invariant. 
	
	
	\begin{definition}
		For $p\in [1,\infty]$, a metric $d$ on $X$ is called a \emph{$p$-metric} if $d(x,y)\leq \elp{\big(d(x,z),d(z,y)\big)}{p}$ for all $x,y,z\in X$.
	\end{definition}
	
	Note that a $1$-metric is just a metric while an $\infty$-metric is an ultrametric (also known as a non-Archimedean metric). Since the $\ell^p$ norms are decreasing in $p$, we see that a $q$-metric is also a $p$-metric whenever $1\leq p\leq q\leq \infty$. 
	By induction, a $p$-metric satisfies $d(x,y) \leq \elp{\left(d(x,z_1),d(z_1,z_2),\dots, d(z_n,y)\right)}{p}$ for any $n\in \N$ and $x,y,z_1,\dots,z_n\in X$.
	
	
	\begin{definition}{\cite[Definition 3.12]{bubenik2019universality}} Let $(X,d,A)$ be a metric pair. For each $p\in [1,\infty]$, we define a new metric $d_p$ on $X$ according to
		\begin{equation}\label{eq:p_strengthened_metric} d_p(x,x') = \min\big(d(x,x'),\elp{\left(d(x,A),d(x',A)\right)}{p}\big),\end{equation}
		for all $x,x'\in X$. We refer to $d_p$ as the \emph{$p$-strengthening of $d$ with respect to $A$}.
	\end{definition}
	
	%
	
	It is an easy observation that $d_p = q^*\overline{d}_p$, where $\overline{d}_p$ is the $p$-quotient metric and $q:X\to X/A$ is the quotient map (see Section \ref{sec:metric-spaces}). Thus $d_p$ is indeed a metric. Moreover, it is clear from the definition that $d_p\leq d$. The following proposition shows that passing from $d$ to $d_p$ leaves the Wasserstein distance unchanged.
	
	\begin{proposition}\label{prop:wasserstein_pstrengthening} For any $p\in [1,\infty]$ we have $W_p[d] = W_p[d_p]$  and $W_p[d](\iota(x),\iota(y)) = d_p(x,y)$, where $\iota:X\to D(X,A)$ is the composition $X\to D(X)\to D(X,A)$ of the inclusion map with the quotient map.
		\begin{proof} 
			To prove the first statement, note that since $d_p\leq d$, it follows immediately from the definition of the Wasserstein distances that $W_p[d_p]\leq W_p[d]$. To prove the reverse inequality, let $\alpha,\beta\in D(X,A)$, let $\epsilon>0$ be given, and let $\sigma = (x_1,x_1') + \dots + (x_r,x_r')\in D(X\times X)$ be a matching between $\alpha$ and $\beta$ with $\Cost_p[d_p](\sigma)<W_p[d_p](\alpha,\beta) + \epsilon$. Define $\sigma'\in D(X\times X)$ as follows. For each $i$, if $d_p(x_i,x_i') \neq d(x_i,x_i')$, so that $d_p(x_i,x_i') = \elp{(d(x_i,A),d(x_i',A))}{p}$, then replace the term $(x_i,x_i')$ in $\sigma$ with $(x_i,a_i) + (a_i',x_i')$, where $a_i,a_i'\in A$ are chosen such that $\elp{(d(x_i,a_i),d(a_i',x_i'))}{p}\leq d_p(x_i,x_i') + \epsilon/r$. Otherwise do nothing. Note that making such a substitution still retains a matching between $\alpha$ and $\beta$. Moreover, we have
			\[\textup{Cost}_p[d](\sigma') \leq \elp{\left(d_p(x_i,x_i') + \epsilon/r\right)_{i = 1}^r}{p} \leq \textup{Cost}_p[d_p](\sigma) + \elp{(\epsilon/r,\dots, \epsilon/r)}{p} \leq \textup{Cost}_p[d_p](\sigma) + \epsilon,\]
			from which the result follows.

			To prove the second statement, note that $(x,y)$ and $(x,a) + (a',y)$ are matchings between $\iota(x)$ and $\iota(y)$ for any choice of $a,a'\in A$. It follows that $W_p[d](\iota(x),\iota(y))\leq d_p(x,y)$. On the other hand, any matching between $\iota(x)$ and $\iota(y)$ contains a term of one of these two forms and hence has cost at least $d_p(x,y)$, proving the reverse inequality.
		\end{proof}
	\end{proposition}
	
	The following is known for the classical Wasserstein distances \cite{mainini2011}. We give a direct proof for persistence diagrams here.
	
	\begin{lemma}[Translation Subinvariance]\label{lem:translation_subinvariance}For any metric pair $(X,d,A)$, $p\in [1,\infty]$, and $\alpha,\beta,\gamma\in D(X,A)$, we have $W_p[d](\alpha + \gamma,\beta + \gamma) \leq W_p[d](\alpha,\beta)$.
		\begin{proof}
			Let $\alpha = x_1 + \dots + x_n$, $\beta = x_1' + \dots + x_m'$, $\gamma = x_1'' + \dots + x_r''\in D(X,A)$. Let $\epsilon>0$ be given. Then there exists a matching $\sigma\in D(X\times X)$ between $\alpha$ and $\beta$ with $\textup{Cost}_p(\sigma)< W_p[d](\alpha,\beta) + \epsilon$. Let $\tilde{\sigma} = \sigma + (x_1'',x_1'') + \dots + (x_r'', x_r'')$. Then $(\pi_1)_*\tilde{\sigma} = \alpha + \gamma \pmod{D(A)}$, $(\pi_2)_*\tilde{\sigma} = \beta + \gamma \pmod{D(A)}$, and
			\[ W_p[d](\alpha + \gamma, \beta + \gamma) \leq \textup{Cost}_p(\tilde{\sigma}) = \textup{Cost}_p(\sigma) < W_p[d](\alpha,\beta) +\epsilon,\]
			which gives the result.
		\end{proof}	
	\end{lemma}
	We can now state the main result of this section. 
	
	\begin{theorem}\label{thm:W_p_ti_when_dp_is_p_metric}
		Let $(X,d,A)$ be a metric pair and let $p\in[1,\infty]$. Then the following are equivalent.
		\begin{enumerate}
			\item $W_p[d]$ is translation invariant;
			\item $d_p:X\times X\to [0, \infty]$ is a $p$-metric.
		\end{enumerate}
		\begin{proof}			
			$(1 \Rightarrow 2)$ Suppose that $W_p[d]$ is translation invariant and let $x,y,z\in X$. Let $\iota:X\to D(X,A)$ denote the composition $X\to D(X)\to D(X,A)$ of the inclusion map and quotient map. Then by Proposition \ref{prop:wasserstein_pstrengthening} and $p$-subadditivity of $W_p[d]$ we have
			\begin{multline*} d_p(x,y) = W_p[d](\iota(x),\iota(y)) = W_p[d](\iota(x) + \iota(z),\iota(z) + \iota(y)) \\\leq \elp{\big(W_p[d](\iota(x),\iota(z)),W_p[d](\iota(z),\iota(y))\big)}{p} = \elp{\big(d_p(x,z),d_p(z,y)\big)}{p},\end{multline*}
			and so $d_p$ is a $p$-metric.
			
			
			$(2\Rightarrow 1)$  
			Assume that $d_p$ is a $p$-metric.
			By \cref{lem:translation_subinvariance}, $W_p[d](\alpha+\gamma,\beta+\gamma) \leq W_p[d](\alpha,\beta)$.
			So it remains to show that 
			$W_p[d](\alpha,\beta) \leq W_p[d](\alpha+\gamma,\beta+\gamma)$.
			By \cref{prop:wasserstein_pstrengthening}, we may instead show that 
			$W_p[d_p](\alpha,\beta) \leq W_p[d_p](\alpha+\gamma,\beta+\gamma)$.
			
			Let $\eps > 0$. Then there is a matching $\sigma \in D(X \times X)$ between $\alpha+\gamma$ and $\beta+\gamma$ such that $\Cost_p[d_p](\sigma) \leq W_p[d_p](\alpha+\gamma,\beta+\gamma) + \eps$.
			Now, consider $(\pi_1)_*\sigma$ and $(\pi_2)_*\sigma$ as sets, where we use indices to avoid higher multiplicities.
			Then $\sigma$ may be viewed as bipartite perfect matching between these two sets. 
			Next, add edges to this bipartite graph between all corresponding (indexed) elements of $\gamma$. 
			We now have a bipartite graph such that
			vertices in $\gamma$ have degree 2 and all other vertices have degree 1.
			Let $n$ be the number of connected components of this graph. 
			Then the connected components of the graph produce a partition of $\sigma$, 
			\begin{equation*}
				\sigma = \sigma_1 + \cdots + \sigma_n, \quad \sigma_i = (x_{i,1},y_{i,1}) + \cdots + (x_{i,m_i},y_{i,m_i}), \ i = 1,\ldots,n,
			\end{equation*}
			where for $i=1,\ldots,n$, $\ x_{i,1} \in \alpha + A$, $y_{i,m_i} \in \beta + A$, and all other terms in $\sigma_i$ lie in $\gamma$.
			Let $\sigma' = (x_{1,1},y_{1,m_1}) + \cdots + (x_{n,1},y_{n,m_n})$.
			Then 
			$(\pi_1)_*\sigma' = \alpha \mod A$ and 
			$(\pi_2)_*\sigma' = \beta \mod A$.
			That is, $\sigma'$ is a matching of $\alpha$ and $\beta$.
			Since 
			$d_p(x_{i,1},y_{i,m_i}) \leq \norm{(d_p(x_{i,1},y_{i,1}),\ldots,d_p(x_{i,m_i},y_{i,m_i}))}_p$ for $i=1,\ldots,n$,
			it follows that 
			\begin{equation*}
				\Cost_p[d_p](\sigma') = \norm{(d_p(x_{i,1},y_{i,m_i}))_{i=1}^n}_p 
				\leq \norm{\left(\norm{(d_p(x_{i,j},y_{i,j}))_{j=1}^{m_i}}_p\right)_{i=1}^n}_p 
				= \Cost_p[d_p](\sigma).
			\end{equation*}
			Therefore $W_p[d_p](\alpha,\beta) \leq W_p[d_p](\alpha+\gamma,\beta+\gamma)+\eps$, which gives the desired result.
		\end{proof}
	\end{theorem}
	
	\begin{corollary}\label{cor:W_1_always_ti}
		For any metric pair $(X,d,A)$, the metric $W_1[d]$ is translation invariant.
	\end{corollary}

	\begin{examples}
		\begin{enumerate}
			\item Let $G = (V,E,\omega)$ be a (possibly directed) weighted graph, where $\omega:E\to [0,\infty]$ defines the edge weights. For a path $\gamma = (e_1,\dots,e_n)$ in $G$ and $p\in [1,\infty]$, the \emph{$p$-cost of $\gamma$} is $C_p(\gamma) = \elp{(\omega(e_i))_{i = 1}^n}{p}$. Define a metric $\rho_p:V\times V\to [0,\infty]$ by setting 
			\[\rho_p(v,v') = \min\left\{C_p(\gamma) \ | \ \gamma \textup{ a path from $v$ to $v'$}\right\},\] 
			if this set is non-empty, and by setting $\rho_p(v,v') = \infty$ otherwise. Then $\rho_p$ is a $p$-metric on $V$. 
			It follows from Theorem \ref{thm:W_p_ti_when_dp_is_p_metric} that $W_p[\rho_p]$ is translation invariant for any choice of $A\subset V$.
			\item The $p$-reflection distance \cite{elchesen2019reflection} is a $p$-metric on the collection of (equivalence classes of) zigzag modules of a fixed length. This can be seen as a special case of the previous example by defining a directed weighted graph whose vertices are equivalence classes of zigzag modules of fixed length $n$ for which there is a directed edge between two vertices if there one of the corresponding zigzag modules can be obtain from the other via a reflection, and with all edge weights set to $1$.
			\item (Metrics on length spaces cannot be $p$-metrics for $p>1$). Let $(X,d)$ be a length space, i.e., for any $x,y\in X$, we have $d(x,y) = \inf\{\ell(p)\ | \ p:I\to X \textup{ a path from $x$ to $y$}\}$, where $\ell(p)$ denotes the length of the path $p$. To avoid trivialities, assume that there exists $x,y\in X$ with $d(x,y)>0$. In a length space, we always have \emph{approximate midpoints}, i.e., for all $x,y\in X$ and $\epsilon>0$ there is a $z\in X$ such that $0\leq d(x,z) - \frac{1}{2}d(x,y) <\epsilon/2$ and $0\leq d(y,z) - \frac{1}{2}d(x,y) <\epsilon/2$. If $d$ is a $p$-metric for some $p\in [1,\infty]$, then for all $\epsilon>0$ we have
			\begin{multline*}
				\quad \quad \quad d(x,y) \leq \elp{\big(d(x,z),d(z,y)\big)}{p} 
				<\elp{\big(d(x,y)/2 +\epsilon/2, d(x,y)/2 + \epsilon/2\big)}{p} \\= \frac{1}{2}(d(x,y) + \epsilon)\elp{(1,1)}{p} = 2^{1/p-1}(d(x,y) + \epsilon).
			\end{multline*}
			Hence $d(x,y) \leq 2^{1/p-1}d(x,y)$ which is only possible if $p\leq 1$. 
			
			\item (Classical persistence diagrams). Consider the pair $(\R^2, \R^2_\geq)$. Note that $\R^2$ is a length space (and in fact a geodesic metric space) whenever the metric $d$ is induced by a $q$-norm, $1\leq q\leq \infty$, and this is still true when we replace $d$ with $d_p$. Thus by the preceding example, $W_p[d]$ is translation invariant if and only if $p = 1$.
		\end{enumerate}
	\end{examples}

\subsection{Grothedieck group completion}

Having described persistence diagrams of metric pairs,
we now apply the the general results of Section \ref{sec:grothendieck_lipschitz_group} to define \emph{virtual persistence diagrams}.

\begin{theorem}\label{thm:virtual_persistence_diagrams}
  Let $(X,d,A)$ be a metric pair, let $p\in [1,\infty]$, and suppose
  that $d_p$ is a $p$-metric. Then $(D(X,A),\allowbreak W_p[d],+)$ is
  a TICL monoid and $W_p[d]$ extends to a metric
  on the Grothendieck group of $D(X,A)$ so as to obtain a TIAL group.
  \begin{proof}
    By Theorem \ref{thm:W_p_ti_when_dp_is_p_metric}, $W_p[d]$ is
    translation invariant and hence by Proposition
    \ref{prop:ti_implies_CL_or_AL}, $(D(X,A),\allowbreak W_p[d],+)$ is
    a TICL monoid. Then by Proposition
    \ref{prop:grothendick_metric_welldefined+ti} and Corollary
    \ref{cor:grothendieck_is_TIAL}, $W_p[d]$ extends to a metric on
    the Grothendieck group of $D(X,A)$ so as to obtain a TIAL group.
  \end{proof}
\end{theorem}

\begin{corollary} \label{cor:vpd-w1}
  Let $(X,d,A)$ be a metric pair.
  Then  $(D(X,A),\allowbreak W_1[d],+)$ is a TICL monoid and $W_1[d]$ extends to a metric 
  on the Grothendieck group of $D(X,A)$ so as to obtain a TIAL group.
\end{corollary}

\begin{definition} \label{def:vpd}
  The corresponding Grothendieck Lipschitz group in \cref{thm:virtual_persistence_diagrams} is denoted $(K(X,A),W_p[d],+)$ and is called the \emph{group of virtual persistence diagrams}. Elements of $K(X,A)$ are called \emph{virtual persistence diagrams}.
\end{definition}

It is straightforward to see that $K(X,A)$ is the quotient $F(X)/F(A)$, where $F(S)$ denotes the free abelian group generated by a set $S$. Hence elements of $K(X,A)$ can be viewed as $\Z$-linear combinations of elements of $X$, with elements of $A$ acting as the identity. By the results of Section \ref{sec:grothendieck_lipschitz_group}, the Wasserstein distances extend to this space in a natural way that is compatible with the group operation, at least when the distance is translation invariant to begin with.
	
	Virtual persistence diagrams appear in \cite{betthauser2019graded} where they arise as the M\"obius inversion of the \emph{$k$th graded rank function}. In that work the authors prove a stability result for the distance $W_1[d]$, where $d$ is the metric obtained from the $1$-norm on $\R^2$.
	
\begin{example}
  Consider the metric pair $(\R^2_{\leq}, d, \Delta)$, where $d$ denotes the Euclidean distance. We can visualize virtual persistence diagrams in this setting
  by graphing points in $\R^2_{\leq}$ together with a multiplicity which is allowed to be negative.
\end{example}



	\section{Wasserstein distance for signed measures}
        \label{sec:measures}
	The Wasserstein distance studied in the present paper is related to a more general distance defined between Radon measures supported on a metric space. This connection has been made precise \cite{divol2019understanding}, with the Wasserstein distance between persistence diagrams being a special case of a ``partial" Wasserstein distance. In this section, we apply the tools developed in the previous sections to the notions of Wasserstein and partial Wasserstein distance. When applied to the monoid of finite Radon measures on a metric space, we obtain a metric of the form introduced in \cite{mainini2011}.
	
	\subsection{The monoid of Radon measures}
	
	Let $(X,d)$ be a metric space and let $\mu$ be a measure on the Borel $\sigma$-algebra of $X$. We say that $\mu$ is \emph{inner regular} if for any open set $U\subset X$,  $\mu(U) = \sup\{\mu(K) \ | \ K\subset X \textup{ compact, } K\subset U\}$;	\emph{outer regular} if for any Borel set $B$, $\mu(B) = \inf\{\mu(U) \ | \ U\subset X \textup{ open, } B\subset U\}$; and \emph{locally finite} if every point of $X$ has a neighborhood $U$ for which $\mu(U)$ is finite. A \emph{Radon measure} is a measure which is inner regular, outer regular, and locally finite. We will denote the set of all Radon measures by $\pmeas(X)$ and the set of all finite Radon measures on $X$ by $\fpmeas(X)$. 
	$\fpmeas(X)$ is a cancellative monoid with the monoid operation taken to be addition of measures, the zero measure being the monoid identity.
	Recall that for a measurable function $f:(X_1,\Sigma_1,\mu)\to (X_2,\Sigma_2)$ from a measure space to a measurable space, the \emph{pushforward measure $f_*\mu$} is a measure on $(X_2,\Sigma_2)$ defined by $(f_*\mu)(E) = \mu(f^{-1}(E))$ for all $E\in \Sigma_2$.

	\subsection{The transportation problem and the Wasserstein distance between measures}
	In this section, we describe the namesake of the Wasserstein distances between persistence diagrams described in previous sections. This is a classical notion of distance between probability measures, or more generally, between measures of equal mass, and is related to the transportation problem for measures. While the similarity between the two notions is evident, the more general notion of partial transportation and the corresponding partial Wasserstein distance, described in Section \ref{sec:partial_wass}, are needed to make the connection precise \cite{divol2019understanding}.
	The classical Wasserstein distance arises from the transportation problem between measures of equal mass supported on a metric space $(X,d)$. 
	
	\begin{definition}[The Transportation Problem] \label{def:transportation_problem}
		Let $(X,d)$ be a metric space and let $\mu,\nu\in\pmeas(X)$ be two positive Radon measures supported on $X$. A \emph{coupling} between $\mu$ and $\nu$ is a Radon measure $\pi$ on $X\times X$ for which $\pi(E\times X) = \mu(E)$ and $\pi(X\times E) = \nu(E)$ for all Borel sets $E\subset X$. The set of all couplings between $\mu$ and $\nu$ is denoted $\Pi(\mu,\nu)$. For $p\in [1,\infty)$, the \emph{$p$-cost} of a coupling is defined by
		\begin{equation}\label{eq:transportation_cost} C_p(\pi) = \int_{X\times X}d(x,y)^pd\pi(x,y).
			\end{equation}
		The \emph{transportation problem} is to find a coupling between $\mu$ and $\nu$ of minimum cost. The \emph{$p$-Wasserstein distance} between $\mu$ and $\nu$ is given by
		\[ \mW_p[d](\mu,\nu) = \left(\inf_{\pi \in \Pi(\mu,\nu)}C_p(\pi)\right)^{1/p},\]
		if $\Pi(\mu,\nu)\neq \emptyset$, and is set to $\infty$ otherwise.
	\end{definition}

	\begin{remark}\label{rem:equal_mass_not_necessary} The $p$-Wasserstein distance is typically viewed as a partial metric in the sense that it is only defined between measures of equal mass. Note that with our definition, the $p$-Wasserstein distance is an (extended) metric on the space of all measures, with the distance between measures of unequal mass always being infinite. Indeed, if $\pi\in \Pi(\mu,\nu)\neq \emptyset$ then $|\mu| = \mu(X) = \pi(X\times X) = \nu(X) = |\nu|$. Hence if $|\mu| \neq |\nu|$ then $\Pi(\mu,\nu) = \emptyset$ and thus $\mW_p[d](\mu,\nu) = \infty$.
	\end{remark}

	The following duality result applies in the $p = 1$ case.
	\begin{theorem}[Kantorovich-Rubinstein Duality {\cite{kellerer1982duality}\cite[Theorem 4.1]{edwards2011KRtheorem}}]\label{thm:KR_duality} Let $\mu,\nu\in\pmeas(X)$ be Radon measures on $(X,d)$ which satisfy the \emph{finiteness conditions} $\int_Xd(x,x_0)d\mu(x), \int_Xd(x,x_0)d\nu(x) <\infty$ for some $x_0\in X$. Then
		\[\textstyle \mW_1[d](\mu,\nu)  = \sup\big\{\int_Xfd\mu - \int_X f d\nu \ \big| \ f:X\to \R, \, \lipnorm{f}\leq 1\big\}.\]
	\end{theorem}

\begin{remark} It is usually assumed in the statement of the previous theorem that $|\mu| = |\nu|$. This is not necessary with our definition of Wasserstein distance. For if $|\mu|>|\nu|$ then, for any $c\in \R$ we have $\int_X c d\mu - \int_X c d\nu = c(|\mu|-|\nu|)\to \infty$ as $c\to\infty$. Since the constant function $c:X\to\R$ is $1$-Lipschitz, the right hand side is infinite. On the other hand, the left hand side is infinite by definition, since $\Pi(\mu,\nu) = \emptyset$ (see Remark \ref{rem:equal_mass_not_necessary}).
	\end{remark}
	
	It follows immediately from Kantorovich-Rubinstein duality that the $1$-Wasserstein distance is translation invariant with respect to addition of measures. This fact together with Proposition \ref{prop:ti_implies_CL_or_AL} immediately implies the following.
	\begin{corollary}\label{cor:measures_are_CL_monoid} Let $(X,d)$ be a metric space. Then	$(\fpmeas(X),\mW_1[d],+)$ is a TICL monoid.
	\end{corollary}
	
	
		In the next section, we will also make use of a generally weaker formulation of the mass transportation problem for $p = 1$, called the \emph{transshipment problem} \cite{edwards2011KRtheorem}\cite[Chapter 6]{rachev1998mass}.

\begin{definition}[The Transshipment Problem] \label{def:transshipment}
  Let $(X,d)$ be a metric space and let $\mu,\nu$ be two Radon measures supported on $X$ with finite first moments.
A \emph{weak coupling} between $\mu$ and $\nu$ is a
  Radon measure $\pi$ on $X\times X$ for which  $\pi(E\times X) +\nu(E) = \pi(X\times E) + \mu(E)$ for all $E\subset X$ Borel.
  Denote by $\Gamma(\mu,\nu)$ the set of all weak couplings between $\mu$ and $\nu$.
  The \emph{$p$-cost} of $\pi$ is defined as for couplings (Definition \ref{def:transportation_problem}).
	\end{definition}

        When the cost is defined by \eqref{eq:transportation_cost} and $p=1$, the transshipment problem and the transportation problem are equivalent \cite[Theorem 4.5]{edwards2011KRtheorem}, i.e.
	\begin{equation}\label{eq:transportation_equiv_transshipment}\inf_{\pi \in \Pi(\mu,\nu)}C_1(\pi) = \inf_{\pi \in \Gamma(\mu,\nu)}C_1(\pi).\end{equation}
	This is no longer true if $d$ is replaced by a more general function, but we will not consider more general cost functions here.
	\subsection{Extending the Wasserstein distance to signed measures} 
	Since the $1$-Wasserstein distance is translation invariant on the monoid $\fpmeas(X)$ of finite Radon measures supported on $(X,d)$, it follows from the results of Section \ref{sec:grothendieck_lipschitz_group} that the $1$-Wasserstein distance extends to the corresponding Grothendieck group. In this section, we describe the Grothendieck group of $\fpmeas(X)$ and describe a transportation formulation of the resulting metric on this group.
	
	The following proposition is essentially the Jordan decomposition theorem.
	\begin{proposition}\label{prop:groth_group_of_measures}
		The Grothendieck group of $\fpmeas(X)$ is the group $\fmeas(X)$ of finite signed Radon measures.
		\begin{proof}
			By the Jordan decomposition theorem, every signed measure $\mu\in \fmeas(X)$ can be written uniquely as a difference $\mu = \mu^+- \mu^-$, where $\mu^+,\mu^-$ are mutually singular positive measures. Define the map $\phi:\fmeas(X)\to K(\fpmeas(X))$ by $\mu = \mu^+- \mu^- \mapsto \mu^+-\mu^-$ (here, the expression on the left hand side denotes elementwise subtraction of positive measures, while the expression on the right hand side represents an equivalence class in $K(\fpmeas(X))$). The verification that $\phi$ is an isomorphism is now straightforward.
		\end{proof}
	\end{proposition}
	
	Since $(\fpmeas(X),\mW_1[d],+)$ is a TICL monoid by Corollary \ref{cor:measures_are_CL_monoid}, by Proposition \ref{prop:uniqueness_of_groth_metric}, $\mW_1[d]$ extends to a unique metric on $\fmeas(X)$.
	\begin{definition}\label{def:signed_wasserstein}
		Let $\mu = \mu^+ - \mu^-$ and $\nu = \nu^+ - \nu^-$ be finite signed measures, written in their Jordan decompositions. The \emph{$1$-Wasserstein distance} between $\mu$ and $\nu$ is defined by
		\[ \mW_1[d](\mu,\nu) = \mW_1[d](\mu^+ + \nu^-, \mu^- + \nu^+).\]
	\end{definition}
	
	Note that the above formula is precisely the global cost extension of the $1$-Wasserstein distance \cite{mainini2011}. Next, we describe a transportation formulation of the Wasserstein distance between signed measures. 
	\begin{definition}[The Signed Transportation Problem]  \label{def:signed-transportation}
		Let $\Delta:X\to X\times X$ be the map $x\mapsto (x,x)$. For $\mu,\nu\in \fmeas(X)$, let  $\Sigma(\mu,\nu)$ denote the collection of all $\sigma\in \fmeas(X\times X)$ of the form  
		\[\sigma = \pi -\Delta_{*}\gamma\]
		for some $\pi\in \fpmeas(X\times X)$ and $\gamma\in \fpmeas(X)$ which satisfy 
		\[\sigma(E\times X) = \mu(E)\quad \textup{ and }\quad\sigma(X\times E) = \nu(E)\]
		for all Borel subsets $E\subset X$. The \emph{cost} of $\sigma\in \Sigma(\mu,\nu)$ is 
		\[ C(\sigma) =  \int_{X\times X}d(x,y)d\sigma(x,y).\]
		The \emph{signed transportation problem} for $\mu,\nu\in \fmeas(X)$ is to find an element of $\Sigma(\mu,\nu)$ of minimum cost.
	\end{definition}
	
	The relationship between the signed transportation problem and the Wasserstein distance between signed measures is the following.
	
	\begin{theorem} \label{thm:signed-transportation}
		For $\mu,\nu\in \fmeas(X)$ we have $\mW_1[d](\mu,\nu) = \inf_{\sigma\in \Sigma(\mu,\nu)} C(\sigma)$.  
		\begin{proof}
			By definition, $\mW_1[d](\mu,\nu) = \inf_{\pi\in \Pi(\mu^+ + \nu^-, \mu^- + \nu ^+)}C_1(\pi)$. Let $\sigma\in \Sigma(\mu,\nu)$ with $\sigma = \pi -\Delta_{*}\gamma$
			for some $\pi\in \fpmeas(X\times X)$ and $\gamma\in \fpmeas(X)$. Then 
			\[ \mu(E)-\nu(E) = \sigma(E\times X)-\sigma(X\times E) = \pi(E\times X)-\pi(X\times E).\]
			Therefore, $\pi\in \Gamma(\mu^++\nu^-,\nu^++\mu^-)$. By change of variables, 
			\[\int_{X\times X}d(x,y)d\Delta_{*}\gamma(x,y) = \int_{X} d(x,y)\Delta(x)d\gamma(x) = \int_{X} d(x,x)d\gamma(x) = 0\]
			and hence we have $C_1(\pi) = \int_{X\times X}d(x,y)d\pi(x,y)  = \int_{X\times X}d(x,y)d\sigma(x,y) = C(\sigma)$. Thus, using the transshipment formulation of the $1$-Wasserstein distance (Equation \eqref{eq:transportation_equiv_transshipment}), we see that 
			\[ \mW_1[d](\mu,\nu) \leq \inf_{\sigma\in \Sigma(\mu,\nu)}C(\sigma).\]
			On the other hand, let $\pi\in \Pi(\mu^++\nu^-,\nu^++\mu^-)$. Then $\sigma =\pi-\Delta_{*}(\mu^-+\nu^-)$ satisfies
			\begin{multline*}\sigma(E\times X) = \pi(E\times X) - \Delta^*(\mu^-+\nu^-)(E\times X) = \pi(E\times X) - \mu^-(E) - \nu^-(E)\\
				 = \mu^+(E) + \nu^{-}(E)- \mu^-(E) - \nu^-(E) = \mu(E),
			\end{multline*}
			\[ \inf_{\sigma\in \Sigma(\mu,\nu)} C(\sigma)\leq   \mW_1[d](\mu,\nu).\qedhere\]
		\end{proof}
	\end{theorem}

	\subsection{The partial transportation problem and partial Wasserstein distances}\label{sec:partial_wass}	
	Let $(X,d,A)$ be a metric pair with $A\subset X$ Borel. The partial transport problem is concerned with measures supported on $X \setminus A$. Informally, the problem is to find the most efficient way to transport the mass of one measure to another, allowing mass to be borrowed from or pushed to the subspace $A$. Expanding upon work of Figalli and Gigli \cite{FIGALLI2010107}, Divol and Lacombe \cite{divol2019understanding} defined the problem for the case that $X$ is some Euclidean space and $A$ is a closed subset of $X$. Motivated by their work, and by adopting the perspective taken in the present paper on the Wasserstein distance between persistence diagrams, we generalize their definition of the partial transport problem and the resulting metrics. By taking $A = \emptyset$, we recover the classical Wasserstein distances as well.
		
	Instead of working with measures on $X \setminus A$ directly, we will instead mimic the algebraic setup we established for persistence diagrams in Section \ref{sec:pers_wass} and consider the quotient $\pmeas(X)/\pmeas(A)$. Here, we are viewing measures $\mu$ supported on $A$ as also being measures on $X$ according to $\mu(E) = \mu(E\cap A)$ for all $E\subset X$ Borel.
	
	\begin{definition}\label{def:relative_measures}
		We define the commutative monoid of \emph{Radon measures on $X$ relative to $A$} by $\pmeas(X,A) = \pmeas(X)/\pmeas(A)$. Similarly, we define the commutative monoids $\fpmeas(X,A) = \fpmeas(X)/\fpmeas(A)$ of \emph{finite Radon measures on $X$ relative to $A$} and $\fmeas(X,A)= \fmeas(X)/\fmeas(A)$ of \emph{finite signed Radon measures on $X$ relative to $A$}.
	\end{definition}
	
	Elements of $\pmeas(X,A)$ are cosets of the form $\mu + \pmeas(A)$, and similarly for $\fmeas(X,A)$ and $\fpmeas(X,A)$. Note that $\mu = \nu \pmod{\pmeas(A)}\iff \mu + \alpha = \nu + \beta$ for some $\alpha,\beta\in \pmeas(A)$, and similarly for $\fmeas(X,A)$ and $\fpmeas(X,A)$.

		\begin{proposition}\label{prop:grothendieck_of_quotient} Let $(X,d,A)$ be a metric pair with $A$ Borel. Then  $K(\fpmeas(X,A))\cong \fmeas(X,A)$.
			\begin{proof}
					By Proposition \ref{prop:groth_group_of_measures}, $K(\fpmeas(X))\cong \fmeas(X)$ and $K(\fpmeas(A)) \cong \fmeas(A)$. Since the monoid of measures is cancellative, it follows from Proposition \ref{prop:groth_of_quotient} that \[K(\fpmeas(X,A))\cong K(\fpmeas(X))/K(\fpmeas(A))\cong \fmeas(X)/\fmeas(A) = \fmeas(X,A).\qedhere\]
			\end{proof}
		\end{proposition}

	\begin{lemma}\label{lem:measure_decomposition}
		Let $A\subset X$ Borel and $\mu\in \pmeas(X)$. Then there exists unique measures $\mu_A\in\pmeas(A)$ and $\mu_C\in \pmeas(X \setminus A)$ such that $\mu = \mu_A + \mu_C$. Analogous statements hold for $\fmeas(X)$ and $\fpmeas(X)$.
		\begin{proof}
			For $E\subset A$ Borel, define $\mu_A(E) = \mu(E)$, and for $E'\subset X \setminus A$ Borel, define $\mu_C(E') = \mu(E')$. View $\mu_A$ as a measure on $X$ by defining $\mu_A(E'') = \mu_A(E''\cap A)$ for $E''\subset X$ Borel, and similarly for $\mu_C$. Then $\mu_A(E'') + \mu_C(E'') = \mu(E''\cap A) + \mu(E''\cap (X \setminus A)) = \mu(E'')$ for all $E''\subset X$ Borel. To see that this decomposition is unique, suppose that $\mu = \mu_A' + \mu_C'$ for some $\mu_A'\in\pmeas(A)$ and $\mu_C'\in \pmeas(X \setminus A)$. For $E\subset A$ Borel, we have $\mu_A(E) = \mu_A(E) + \mu_C(E) = \mu(E) = \mu_A'(E) + \mu_C'(E) = \mu_A'(E)$, since $\mu_C(E) = \mu_C'(E) = 0$, and hence $\mu_A = \mu_A'$. A similar argument shows that $\mu_C = \mu_C'$.The proofs for $\fmeas(X)$ and $\fpmeas(X)$ are identical.
		\end{proof}
	\end{lemma}
	\begin{proposition}\label{prop:quotient_of_measues_is_measures_on_difference}
		Let $(X,d,A)$ be a metric pair with $A\subset X$ Borel. There are monoid isomorphisms $\pmeas(X,A) \cong \pmeas(X \setminus A)$, $\fpmeas(X,A) \cong \fpmeas(X \setminus A)$, and $\fmeas(X,A) \cong \fmeas(X \setminus A)$.
		\begin{proof}
			Define $\Phi:\pmeas(X \setminus A)\to\pmeas(X)/\pmeas(A)$ by $\nu\mapsto \nu + \pmeas(A)$. Clearly $\Phi$ is a monoid homomorphism. To see that $\Phi$ is surjective, let $\mu\in \pmeas(X)$ be given. By Lemma \ref{lem:measure_decomposition}, there exists unique measures $\mu_A\in\pmeas(A)$ and $\mu_C\in \pmeas(X \setminus A)$ such that $\mu = \mu_A + \mu_C$. Thus $\mu = \mu_C \pmod{\pmeas(A)}$ and hence $\Phi(\mu_C) = \mu + \pmeas(A)$. To see that $\Phi$ is injective, suppose that $\Phi(\nu) = \Phi(\nu')$. Then there are measures $\alpha,\beta\in \pmeas(A)$ such that $\nu + \alpha = \nu' + \beta$. Since $\nu,\nu'\in \pmeas(X \setminus A)$ and $\alpha,\beta\in \pmeas(A)$, by the uniqueness statement of Lemma \ref{lem:measure_decomposition} we have $\nu = \nu'$, as desired. The second isomorphism is proven similarly. The last isomorphism then follows from the second and Propositions \ref{prop:groth_group_of_measures} and \ref{prop:grothendieck_of_quotient}.
		\end{proof}
	\end{proposition}

We now introduce a generalization of the Wasserstein distances, of which the classical Wasserstein distances can be seen as a special case.
	
	\begin{definition} \label{def:wasserstein-radon-pairs}
		Let $p\in [1,\infty)$ and let $\mu,\nu\in \pmeas(X)$. We say that $\mu,\nu$ are \emph{$p$-finite with respect to $A$} if
		\begin{equation}\label{eq:finiteness} \int_{X}d(x,A)^pd\mu(x) <\infty,\qquad \int_{X}d(x,A)^pd\nu(x)< \infty.
		\end{equation}
		Let $\mathcal{M}_{p,A}^+(X)$ denote the submonoid of $\pmeas(X)$ consisting of all measures that are $p$-finite with respect to $A$.
Let $\mathcal{M}_p^+(X,A) = \mathcal{M}_{p,A}^+(X)/\pmeas(A)$.
                  For $\mu,\nu\in \pmeas_{p,A}(X)$, we denote by $\Pi_A(\mu,\nu)$ the subset of $\pmeas(X\times X)$ of all Radon measures $\sigma$ supported on $X\times X$ which satisfy
		\begin{equation*}
			(\pi_1)_*\sigma = \mu \pmod{\pmeas(A)}, \qquad (\pi_2)_*\sigma = \nu \pmod{\pmeas(A)}.
		\end{equation*}
		The \emph{$p$-cost} of $\sigma\in \Pi_A(\mu,\nu)$, denoted $C_p^A(\sigma)$, is defined as in the ordinary transportation cost \eqref{eq:transportation_cost}. The \emph{partial optimal transportation distance between $\mu$ and $\nu$ with respect to $A$} is then defined by
		\[ \mW_p^A[d](\mu,\nu) = \left(\inf_{\sigma\in \Pi_A(\mu,\nu)} C_p^A(\sigma)\right)^{1/p}\]
		when $\Pi_A(\mu,\nu)\neq \emptyset$, and we set $W^A_p(\mu,\nu) = \infty$ otherwise.
	\end{definition}

\begin{remarks} \label{rem:wasserstein-radon-pairs}
		\begin{enumerate}
			\item The $p$-finite conditions are analogous to the finiteness conditions in Theorem \ref{thm:KR_duality} and guarantee that the partial optimal transport distance is finite between measures satisfying this condition.
			\item $C_p^A$ can also be defined as an integral over $(X\times X)\setminus(A\times A)$ without changing the Wasserstein distance, since measures are only considered modulo measures on $A$.
			\item $\mW_p^A[d]$ is zero on $\pmeas(A)$ and thus induces a distance on $\pmeas_p(X,A)$ which we denote $\mW_p[d]$. By Proposition \ref{prop:quotient_of_measues_is_measures_on_difference}, we can alternatively view $\mW_p[d]$ as a metric between Radon measures supported on $X\setminus A$ that are $p$-finite with respect to $A$.
		\end{enumerate}
	\end{remarks}

	Divol and Lacombe define the partial optimal transport distance slightly differently for measures supported on an open subset $\mathcal{S}\subset \R^d$. This distance was originally introduced by Figalli and Gigli for bounded subsets of Euclidean spaces \cite{FIGALLI2010107}.
	
	\begin{definition}[Divol and Lacombe, Figalli and Gigli \cite{divol2019understanding,FIGALLI2010107}]\label{def:partial_wass_divol_lacombe}
		Let $p\in[1,\infty)$, let $\mathcal{S}\subset \R^d$ be open, and let $\mu,\nu$ be Radon measures supported on $\mathcal{S}$ which satisfy
		\begin{equation}\label{eq:dl_finiteness}
			\int_\mathcal{S}d(x,\partial \mathcal{S})^p d\mu(x)<\infty,\quad \int_{\mathcal{S}}d(x,\partial \mathcal{S})^p d\nu(x)<\infty.
		\end{equation}
		Let $\textup{Adm}(\mu,\nu)$ denote the set of all Radon measures $\sigma$ on $\overline{\mathcal{S}}\times\overline{\mathcal{S}}$ which satisfy
		\[\sigma(E\times\overline{\mathcal{S}}) = \mu(E),\qquad \sigma(\overline{\mathcal{S}}\times E') = \nu(E'),\]
		for all $E,E'\subset \mathcal{S}$ Borel. The \emph{$p$-cost} of $\sigma\in \textup{Adm}(\mu,\nu)$ is defined by
		\[ \textup{C}_p^{\textup{OT}}(\sigma) = \int_{\overline{\mathcal{S}}\times\overline{\mathcal{S}}}d(x,y)^pd\sigma(x,y).\]
		Then define 
		\[ \textup{OT}_p(\mu,\nu) =  \left(\inf_{\sigma\in \textup{Adm}(\mu,\nu)} C_p^{\textup{OT}}(\sigma)\right)^{1/p}.\]
	\end{definition}
	
	\begin{lemma}\label{lem:distance_to_boundary}
		For $A\subset \R^d$ and $x\not\in A$ we have $d(x,A) = d(x,\partial A)$.
		\begin{proof}
			Since $\overline{A} = \textup{int}(A)\cup \partial A$ is closed and $d(x,\overline{A}) = d(x,A)$, there exists $a\in \overline{A}$ such that $d(x,A) = d(x,a)$. Suppose that $a\in \textup{int}(A)$. Let $\epsilon>0$ be such that $\overline{B}_\epsilon(a)\subset \textup{int}(A)$ and let $a'= a +
			\epsilon\frac{x-a}{\|x-a\|}$. Then $a'\in \overline{B}_\epsilon(a)\subset A$ and $d(x,a') = d(x,a) -\epsilon <d(x,A)$, a contradiction. 
		\end{proof}
	\end{lemma}
	Recall that a \emph{retraction} of a topological space $X$ onto a subspace $A\subset X$ is a continuous map $r:X\to A$ such that $r(a) = a$ for all $a\in A$.
	\begin{proposition}\label{prop:OT_equals_partial} 
		Let $\mathcal{S}\subset \R^d$ be an open subset and let $\mu,\nu$ be Radon measures supported on $\mathcal{S}$ which satisfy \eqref{eq:dl_finiteness}. Let $X = \R^d$ and $A = X\setminus \mathcal{S}$. Suppose also that there is a $1$-Lipschitz retraction $r:X\to \overline{\mathcal{S}}$ with $r(A)\subset \partial A$. Then
		\[\mW_p^A[d](\mu,\nu) = \textup{OT}_p(\mu,\nu),\]
		where on the left-hand-side, we view $\mu$ and $\nu$ as measures supported on $X$ in the usual way.
		\begin{proof}
			By Lemma \ref{lem:distance_to_boundary} and the fact that $d(x,A) = 0$ for all $x\in A$, we have
			\[ \int_X d(x,A)^pd\mu(x) = \int_\mathcal{S}d(x,A)^pd\mu(x) = \int_\mathcal{S}d(x,\partial A)^pd\mu(x) = \int_{\mathcal{S}}d(x,\partial\mathcal{S})^pd\mu(x),\]
			and similarly for $\nu$. Thus $\mu,\nu$ satisfy \eqref{eq:finiteness} if and only if they satisfy \eqref{eq:dl_finiteness}. Now let $\sigma\in \textup{Adm}(\mu,\nu)$. View $\sigma$ as an element of $\pmeas(X\times X)$ by setting $\sigma(U) = \sigma(U\cap (\overline{\mathcal{S}}\times \overline{\mathcal{S}}))$ for all $U\subset X\times X$ Borel. Then, for the projection $\pi_1:X\times X\to X$, we have 
			\begin{multline*}
				(\pi_1)_*\sigma(E) = \sigma(E\times X) = \sigma((E\times X)\cap(\overline{\mathcal{S}}\times \overline{\mathcal{S}})) 
				 = \sigma((E\cap\overline{\mathcal{S}})\times \overline{\mathcal{S}})\\ = \sigma((E\cap\mathcal{S})\times \overline{\mathcal{S}}) + \sigma((E\cap\partial \mathcal{S})\times \overline{\mathcal{S}})
				 =  \mu(E\cap\mathcal{S}) + \sigma((E\cap\partial \mathcal{S})\times \overline{\mathcal{S}})
		\end{multline*}
		 for $E\subset X$ Borel. It follows that $(\pi_1)_*\sigma(E) + \mu(E\cap A) = \mu(E\cap\mathcal{S}) + \mu(E\cap A) + \sigma((E\cap\partial \mathcal{S})\times \overline{\mathcal{S}}) = \mu(E) + \sigma((E\cap\partial \mathcal{S})\times \overline{\mathcal{S}})$. Now let $\alpha,\beta\in\pmeas(A)$ be defined by $\alpha(F) := \mu(F\cap A)$ and $\beta(F) := \sigma((F\cap \partial \mathcal{S})\times\overline{\mathcal{S}})$ for all $F\subset A$ Borel. View $\alpha$ as a measure on $X$ by defining $\alpha(E') = \alpha(E'\cap A)$ for $E'
		 \subset X$ Borel, and similarly for $\beta$. Then we have 
		\[(\pi_1)\sigma_*(E) + \alpha(E) = (\pi_1)_*\sigma(E) + \mu(E\cap A) = \mu(E) + \sigma((E\cap\partial\mathcal{S})\times\overline{\mathcal{S}}) = \mu(E) + \beta(E),\]
		 for all $E\subset X$ Borel. Thus $(\pi_1)_*\sigma = \mu \pmod{\pmeas(A)}$. Similarly, we have $(\pi_2)_*\sigma = \nu \pmod{\pmeas(A)}$. Now, since $\sigma$ is supported on $\overline{\mathcal{S}}\times\overline{\mathcal{S}}$, we have
			$\int_{X\times X}d(x,y)^pd\sigma(x,y) = \int_{\overline{\mathcal{S}}\times\overline{\mathcal{S}}}d(x,y)^pd\sigma(x,y)$ so that $C_p(\sigma) = C_p^{\textup{OT}}(\sigma)$. It follows that $\mW_p^A[d](\mu,\nu) \leq \textup{OT}_p(\mu,\nu)$.
			
			To prove the reverse inequality, let $\sigma\in\Pi(\mu,\nu)$ be given. Define $\sigma' = (r\times r)_*\sigma$. Then $\sigma'$ is supported on $\overline{\mathcal{S}}\times\overline{\mathcal{S}}$. To see that $\sigma'\in \textup{Adm}(\mu,\nu)$, let $E\subset \mathcal{S}$ be a Borel set. Then $\sigma'(E\times \overline{\mathcal{S}}) = \sigma((r\times r)^{-1}(E\times \overline{\mathcal{S}})) = \sigma(r^{-1}(E)\times r^{-1}(\overline{\mathcal{S}}))$. Since $r(A)\subset \partial A = \partial\mathcal{S}$ and $\mathcal{S}\subset \mathbb{R}^d$ is open by assumption, we have $r^{-1}(E) = E$. Thus, $\sigma'(E\times\overline{\mathcal{S}}) = \sigma(r^{-1}(E)\times r^{-1}(\overline{\mathcal{S}})) = \sigma(E\times X)$. Since $\sigma\in \Pi(\mu,\nu)$, there are $\alpha,\beta\in \pmeas(A)$ such that $\sigma(E\times X) + \alpha(E)= \mu(E) + \beta(E)$. Since $E\subset\mathcal{S}$, $\alpha(E) = \beta(E) = 0$ so that $\sigma'(E\times \overline{\mathcal{S}}) = \mu(E)$. Similarly,  $\sigma'(\overline{\mathcal{S}}\times E) = \nu(E)$ and thus $\sigma'\in \textup{Adm}(\mu,\nu)$.
			Now, by the change of variables formula together with the fact that $r$ is $1$-Lipschitz, we have
			\[C_p^\textup{OT}(\sigma') = \int_{\overline{\mathcal{S}}\times\overline{\mathcal{S}}}d(x,y)^pd\sigma'(x,y) = \int_{X\times X}d(r(x),r(y))^pd\sigma(x,y) \leq \int_{X\times X}d(x,y)^pd\sigma(x,y) = C_p(\sigma).\]
			Thus $\textup{OT}_p(\mu,\nu)\leq \mW_p^A[d](\mu,\nu)$.
		\end{proof}
		
		\begin{corollary} \label{cor:agrees-dl}
			$\mW_p^A[d] = \textup{OT}_p$ for $\mathcal{S} = \R^2_<\subset \R^2$. 
			\begin{proof}
				Let $X = \R^2$ and $A = X\setminus \mathcal{S} = R^2_\geq$. The map $r:X\to \overline{\mathcal{S}} = \R^2_\leq$ given by $r(\mathbf{x}) = \mathbf{x}$ for $\mathbf{x}\in \R^2_\leq$ and $r((x,y)) = (\frac{x+y}{2},\frac{x+y}{2})$ for $(x,y)\in \R^2_>$ is a $1$-Lipschitz retraction of $X$ onto $\overline{\mathcal{S}}$ with $r(A)\subset \partial A = \Delta$. The result then follows from Proposition \ref{prop:OT_equals_partial}.
			\end{proof}
		\end{corollary}
	\end{proposition}

	%
	
	The classical transportation problem can be seen as the special case of the partial transportation problem in which the subset $A$ is taken to be empty. The Wasserstein distance is then just the partial optimal transportation distance for the metric pair $(X,d,\emptyset)$. In this case, $\textup{Adm}(\mu,\nu)$ reduces to the set of all measures $\pi$ supported on $X\times X$ which satisfy $\pi(E\times X)  = \mu(E)$ and $\pi(X\times E)  = \nu(E)$ for all Borel sets $E\subset X$, i.e. $\textup{Adm}(\mu,\nu) = \Pi(\mu,\nu)$. 
	
	In the case that $X$ is $\R^2$ equipped with a Euclidean distance and $A = \R^2_\geq$, Divol and Lacombe have observed conversely that the partial optimal transport distance can be formulated in terms of the classical Wasserstein distance. 
	
	\begin{proposition}[{\cite[Proposition 3.7]{divol2019understanding}}]\label{prop:partial_wass_from_classical_wass}
          Let $\mu,\nu$ be finite Radon measures supported on $\R^2_<$
          which satisfy \eqref{eq:dl_finiteness}. Let $t\geq \mu(\R^2_<) + \nu(\R^2_<)$. Define measures on the quotient space $\R^2/A$ by $\tilde{\mu} = \mu + (t-\mu(\R^2_<))\delta_{A}$ and $\tilde{\nu} = \nu + (t-\nu(\R^2_<))\delta_{A}$ (here, we are viewing $\mu$ and $\nu$ as measures on the quotient space). Then $\textup{OT}_p(\mu,\nu) = \mW_p[d](\tilde{\mu},\tilde{\nu})$.
	\end{proposition}
	
	The preceding proposition together with Kantorovich-Rubinstein duality imply the following.
	\begin{corollary}\label{cor:OT_ti}
		$\textup{OT}_1$, and hence $\mW_1^{\R^2_\geq}[d]$, is translation invariant on $\pmeas_{1,\R_{\geq}^2}(\R^2_<)$.
		\begin{proof}
			Let $\mu,\nu,\kappa$ be measures supported on $\R^2_<$ which satisfy \eqref{eq:dl_finiteness}. Let $r\geq 2(\mu(\R^2_<) + \nu(\R^2_<) + \kappa(\R^2_<))$ and let $\tilde{\mu} = \mu+ (r/2-\mu(\R^2_<))\delta_A$, $\tilde{\nu} = \nu + (r/2-\nu(\R^2_<))\delta_A$, and $\tilde{\kappa} = \kappa + (r/2 - \kappa(\R^2_<))$. Then by Proposition \ref{prop:partial_wass_from_classical_wass}, $\textup{OT}_1(\mu + \kappa, \nu + \kappa) = \mW_1[d](\tilde{\mu} + \tilde{\kappa}, \tilde{\nu} + \tilde{\kappa})$. By Kantorovich-Rubinstein duality, $\mW_1[d](\tilde{\mu} + \tilde{\kappa}, \tilde{\nu} + \tilde{\kappa}) = \mW_1[d](\tilde{\mu}, \tilde{\nu})$ and again by Proposition \ref{prop:partial_wass_from_classical_wass},  $\mW_1[d](\tilde{\mu},\tilde{\nu}) = \textup{OT}_1(\mu,\nu)$	as desired.
		\end{proof}
	\end{corollary}
	
	\begin{corollary}\label{cor:extending_classical_persistence_wasserstein}
		$\mW_1^{\R^2_\geq}[d]$ extends to a unique metric $\mW_1[d]$ on $\mathcal{M}_1(\R^2,\R^2_\geq) = K(\pmeas_1(\R^2,\R^2_\geq))$.
		\begin{proof}
	By Corollary \ref{cor:OT_ti}, $\mW_1^{\R^2_\geq}[d]$ is translation invariant on $\pmeas_1(\R^2,\R^2_\geq)$. Hence, by Proposition \ref{prop:ti_implies_CL_or_AL}, $(\pmeas_1(\R^2,\R^2_\geq),\mW_1^{\R^2_\geq}[d],+)$ is a TICL monoid and result follows from Proposition \ref{prop:uniqueness_of_groth_metric}.
	\end{proof}
	\end{corollary}
	
\subsection{Persistence diagrams as Radon measures}
\label{sec:pd-radon}

Let $(X,d,A)$ be a metric pair. Recall that there is a canonical embedding from $X \to D(X)$ given by sending $x \in X$ to the indicator function on $X$. Similarly, there is a canonical embedding from $X \to \fpmeas(X)$ given by sending $x \in X$ to the Dirac measure on $x$. This embedding extends by linearity to a canonical embedding $D(X) \to \fpmeas(X)$.
  This embedding induces a canonical embedding $D(X,A) \to \fpmeas(X,A)$.
  That is, persistence diagrams may be viewed as special cases of Radon measures.

\begin{proposition} \label{prop:embedding}
 Let $\alpha,\beta \in D(X)$. Under the canonical embedding $D(X) \to \fpmeas(X)$, we may view $\alpha$ and $\beta$ as measures.
  Then $\mW_p[d](\alpha,\beta) = W_p[d](\alpha,\beta)$. Therefore, we have a canonical isometric embedding $(D(X,A),W_p[d]) \to (\fpmeas(X,A),\mW_p[d])$.
\end{proposition}

\begin{proof}
    By the fundamental theorem of linear programming there is a solution to the transportation problem for $\alpha$ and $\beta$ (viewed as measures) given by a matching of $\alpha$ and $\beta$ (viewed as persistence diagrams).
\end{proof}
  

	

	\section{Infinite persistence diagrams: the Cauchy completion of $D(X,A)$}\label{sec:cauchy_completion_of_pers_diags}
	As we saw in Section \ref{sec:pers_wass}, the space of finite persistence diagrams has a nice description in terms of the free commutative monoid. While persistence diagrams arising in practice are finite, these are sometimes finite approximations of infinite persistence diagrams. Following the works of \cite{mileyko2011probability} and \cite{blumberg2014robust}, in this section we describe a family of infinite persistence diagrams satisfying a certain finiteness condition in terms of the Cauchy completion of the Lipschitz monoid of finite persistence diagrams. The main contribution of the work in this section is threefold: we generalize the results of those papers to pairs of metric spaces, show that the monoid operation of the space of finite diagrams extends to the Cauchy completion, and state a corresponding universal property for the resulting space.
	
	\subsection{Cauchy completion of spaces of diagrams}
	Let $(X,d,A)$ be a metric pair, let $p\in [1,\infty]$, and consider the CL monoid $(D(X,A),W_p[d],+)$. We will describe the Cauchy completion of the metric space $(D(X,A),W_p[d])$ and show that the above monoid operation extends to this completion so as to again obtain a CL monoid. We will first introduce a definition.
	

		\begin{definition}\label{def:quotient_by_dist_zero}
			For a metric space $(X,d)$, let $\sim_0$ denote the equivalence relation on $X$ given by $x\sim_0 y\iff d(x,y) = 0$. The corresponding quotient space is denoted $(X/\!\!\sim_0,\tilde{d})$.
		\end{definition}
		If $(M,d,+)$ is a CL monoid, then addition descends to the quotient $(M/\!\!\sim_0,\tilde{d})$ so as to obtain a CL monoid. Indeed, define $[x] \tilde{+} [y] = [x+y]$. This is well-defined since if $[x] = [x']$ and $[y] = [y']$ then $d(x+y,x'+y') \leq \lipnorm{+}(d(x,x') + d(y,y')) = 0$. Moreover, it is easily verified that $\lipnorm{\tilde{+}} = \lipnorm{+}$.
	
	\begin{definition}
		Let $(X,d)$ be a metric space. Let $\tilde{X}$ denote the set of Cauchy sequences in $X$. Define an equivalence relation $\sim$ on $\tilde{X}$ by declaring that $(x_n)\sim (y_n)\iff \lim_{n\to \infty}d(x_n,y_n)  = 0$. Let $\overline{X} = \tilde{X}/\!\!\sim$ and let $[x_n]$ denote the equivalence class of $(x_n)$ under $\sim$. We define a metric $\overline{d}$ on $\overline{X}$ by setting $\overline{d}([x_n],[x_n']) = \lim_{n\to \infty}d(x_n,x_n')$. The metric space $(\overline{X},\overline{d})$ is the \emph{Cauchy completion} of $(X,d)$. The Cauchy completion of a metric pair $(X,d,A)$ is defined to be $(\overline{X},\overline{d},A/\!\!\sim_0)$.
	\end{definition}
        The Cauchy completion $(\overline{X},\overline{d})$ is complete and contains a dense isometric copy of $(X/\!\!\sim_0,\tilde{d})$, the image of the map $i:X\to \overline{X}$ which sends $x\in X$ to the equivalence class of the constant sequence at $x$. If $X$
        satisfies the separation axiom
        then $i$ is an isometric embedding of $X$ into $\overline{X}$.
	
	We record the following lemma for later reference.
	
	\begin{lemma}\label{lem:cauchy_inclusion_is_epi}
		The canonical map $i:X\to \overline{X}$ is an epimorphism in $\lip$.
		\begin{proof}
			Let $f,g:\overline{X}\to Z$ be Lipschitz and hence continuous. Since $i(X)$ is dense, $f$ and $g$ are completely determined by their values on $i(X)$. Hence, if $f\circ i = g\circ i$ then $f = g$ and thus $i$ is an epimorphism.
		\end{proof}
	\end{lemma} 
	We note, however, that $i$ is a monomorphism if and only if $X$ satisfies the separation axiom. When $X$ does satisfy the separation axiom, $i$ is an example of a morphism which is both monic and epic but is not necessarily an isomorphism unless $X$ is complete to begin with. The following lemma shows that elements of the Cauchy completion can be obtained as limits of constant sequences approximating them.
\begin{lemma}\label{lem:limit_in_cauchy_completion}
			Let $(X,d)$ be a metric space with Cauchy completion $(\overline{X},\overline{d})$ and let $i:X\to \overline{X}$ be the canonical map. For an element $[x_n]\in \overline{X}$ we have $[x_n] = \lim_{m\to\infty}i(x_m)$.
		\begin{proof}
			Given $[x_n]\in\overline{X}$, we have $\lim_{m\to\infty}\overline{d}(i(x_m),[x_n]) = \lim_{m\to\infty}\lim_{n\to\infty}d(x_m,x_n) = 0$, since $(x_n)$ is a Cauchy sequence, and hence $\lim_{m\to\infty}i(x_m) = [x_n]$.
		\end{proof}
		\end{lemma}

	The Cauchy completion, together with the canonical inclusion $i:X\to \overline{X}$, satisfies the following ``enriched" universal property.
	
	\begin{theorem}\label{thm:univ_prop_of_cauchy_completion} Let $(X,d)$ be a metric space. Let $\mathbf{Lip_{cpl,sep}}$ denote the full subcategory of $\mathbf{Lip}$ whose objects are complete metric spaces satisfying the separation axiom. Then 
		\begin{thmenum}
			\item \label{thm:univ_prop_of_cauchy_completion_part_1}For any $(N,\rho)\in \mathbf{Lip_{cpl,sep}}$ and Lipschitz map $\phi:X\to N$, there is a unique Lipschitz map $\overline{\phi}:\overline{X}\to N$ such that $\overline{\phi}\circ i = \phi$; and
			\item \label{thm:univ_prop_of_cauchy_completion_part_2} $\lipnorm{\overline{\phi}} = \lipnorm{\phi}$.
		\end{thmenum}
		\begin{proof}
			(1) Let $\phi:X\to N$ be a given Lipschitz map and define $\overline{\phi}:\overline{X}\to N$ by $\overline{\phi}([x_n])= \lim_{n\to\infty}\phi(x_n)$. Since $\phi$ is Lipschitz, Cauchy sequences are sent to Cauchy sequences. Since $N$ is complete, the limit defining $\overline{\phi}$ exists and is unique since $N$ satisfies the separation axiom. Moreover, $\overline{\phi}(i(x)) = \lim_{n\to \infty}\phi(x) = \phi(x)$ for all $x\in X$ so that $\overline{\phi}\circ i = \phi$. Suppose $\psi:\overline{X}\to N$ is Lipschitz and $\psi\circ i = \phi$. Then by Lemma \ref{lem:limit_in_cauchy_completion}, we have
			\[ \psi([x_n]) = \psi(\lim_{m\to\infty}i(x_m)) = \lim_{m\to\infty}\psi(i(x_m)) = \lim_{m\to\infty}\phi(x_m),\]
			so that $\psi = \overline{\phi}$.
			
			(2) We have
			\begin{multline*}\rho(\overline{\phi}([x_n]),\overline{\phi}([x_n'])) = \rho(\lim_{n\to \infty}\phi(x_n),\lim_{n\to \infty}\phi(x_n')) = \lim_{n\to \infty}\rho(\phi(x_n),\phi(x_n'))\\
				\leq \lipnorm{\phi}\lim_{n\to \infty}d(x_n,x_n') = \lipnorm{\phi}\overline{d}([x_n],[x_n']),
			\end{multline*}
			so that $\overline{\phi}$ is Lipschitz with $\lipnorm{\overline{\phi}}\leq \lipnorm{\phi}$. On the other hand, $\lipnorm{\phi} = \lipnorm{\overline{\phi}\circ i}\leq \lipnorm{\overline{\phi}}\lipnorm{i} = \lipnorm{\overline{\phi}}$, since $\lipnorm{i} = 1$, and thus $\lipnorm{\overline{\phi}} = \lipnorm{\phi}$.
		\end{proof}
	\end{theorem}

The analogous result holds for metric pairs $(X,d,A)$ and its Cauchy completion $(\overline{X},\overline{d},A/\!\!\sim_0)$, with Lipschitz maps replaced with Lipschitz maps of pairs.
	
	\begin{corollary}\label{cor:cauchy_completion_unique_upto_isom}
		If $(\hat{X},\hat{d})$ together with a Lipschitz map $u:X\to\hat{X}$ satisfy both conclusions of Theorem \ref{thm:univ_prop_of_cauchy_completion}, then $(\hat{X},\hat{d})$ and $(\overline{X},\overline{d})$ are isometric.
		\begin{proof}
			Since $i:X\to \overline{X}$ is Lipschitz, there exists a unique Lipschitz map $\hat{i}:\hat{X}\to \overline{X}$ such that $\hat{i}\circ u = i$, and $\lipnorm{\hat{i}} = \lipnorm{i} = 1$. Similarly, there is a unique Lipschitz map $\overline{u}:\overline{X}\to\hat{X}$ such that $\overline{u}\circ i = u$, with $\lipnorm{\overline{u}} = \lipnorm{u}$. Since $\textup{id}_{\hat{X}}$ is the unique Lipschitz map such that $\textup{id}_{\hat{X}} \circ u = u$, we have $\lipnorm{u} = \lipnorm{\textup{id}_{\hat{X}}} = 1$. Thus $\lipnorm{\overline{u}} = \lipnorm{u} = \lipnorm{\textup{id}_{\hat{X}}} = 1$.
			
			Now, again by uniqueness, $\hat{i}\circ \overline{u} = \textup{id}_{\overline{X}}$ and $\overline{u}\circ\hat{i}= \textup{id}_{\hat{X}}$. Since $\hat{i}$ and $\overline{u}$ are both $1$-Lipschitz, they give the desired isometry.
		\end{proof}
	\end{corollary}
	
	\begin{corollary}\label{cor:completion_of_prod_is_prod_of_completion}
		The Cauchy completion of $(X\times X,d+d)$ is (up to isometry) $(\overline{X}\times \overline{X},\overline{d}+\overline{d})$.
		\begin{proof}
			Let $i:X\to \overline{X}$ be the canonical inclusion. Then $i\times i:X\times X\to \overline{X}\times \overline{X}$ is easily seen to satisfy both conclusions of Theorem \ref{thm:univ_prop_of_cauchy_completion}. Hence by Corollary \ref{cor:cauchy_completion_unique_upto_isom}, $(\overline{X}\times\overline{X},\overline{d}+\overline{d})$ is (up to isometry) the Cauchy completion of $(X\times X,d+d)$.
		\end{proof}
	\end{corollary}
	
	We can now see formally that addition on a Lipschitz monoid extends to its Cauchy completion.
	
	\begin{proposition} \label{prop:cauchy_completion_is_cl}
		Let $(M,d,+)$ be a CL monoid. Then the Cauchy completion $(\overline{M},\overline{d})$ admits a CL monoid structure which restricts to that of $M/\!\!\sim_0$ on the isometric embedding of $M/\!\!\sim_0$ in $\overline{M}$. Moreover, the canonical inclusion $i:M\to\overline{M}$ is a monoid homomorphism, and if $d$ is translation invariant then so is $\overline{d}$. 
		\begin{proof}
			Since $M$ is a CL monoid, the map $+:M\times M\to M$ is Lipschitz. By Corollary \ref{cor:completion_of_prod_is_prod_of_completion}, the Cauchy completion of $(M\times M, d+d)$ is $(\overline{M}\times \overline{M},\overline{d}+\overline{d})$, and by the universal property of the Cauchy completion, there is a unique Lipschitz map $+:\overline{M}\times\overline{M} \to \overline{M}$ such that $+\circ (i\times i) = i\circ +$ (here, we denote the operations on $M$ and $\overline{M}$ both by $+$). The fact that $+$ makes $(\overline{M},\overline{d},+)$ a CL monoid is easily verified. Moreover, the equality $+\circ (i\times i) = i\circ +$ means exactly that $i(x) + i(y) = i(x +y)$ and so $i$ is an embedding of $M/\!\!\sim_0$ as a submonoid of $\overline{M}$. We can then compute addition in $\overline{M}$ explicitly as follows. For $[x_n],[y_n]\in \overline{M}$ we have $[x_n]+[y_n] = \lim_n i(x_n) + \lim_n i(y_n) = \lim_n i(x_n) + i(y_n) = \lim_n i(x_n+y_n) = [x_n+y_n]$. Now to see that $\overline{d}$ is translation invariant on $\overline{M}$, if $[x_n],[y_n],[z_n]\in \overline{X}$ then we have $\overline{d}([x_n]+[z_n],[y_n]+[z_n]) = \lim_{n\to \infty}d(x_n+z_n,y_n+z_n) = \lim_{n\to\infty} d(x_n,y_n) = d([x_n],[y_n])$.
		\end{proof}
	\end{proposition}
	The next proposition shows that the unique Lipschitz map provided by universality of the Cauchy completion is a monoid homomorphism when the given Lipschitz map is a monoid homomorphism between CL monoids.
	\begin{proposition}\label{prop:unique_morphism_is_monoid_homomorphism}
		Let $(M,d,+)$ be a given CL monoid. Suppose that $(N,\rho,+)$ is any other CL monoid and that $\phi:M\to N$ is a Lipschitz monoid homomorphism. Then the unique Lipschitz map $\overline{\phi}:\overline{M}\to N$ provided by Theorem \ref{thm:univ_prop_of_cauchy_completion} is a monoid homomorphism when $\overline{M}$ is equipped with the CL monoid structure provided by Proposition \ref{prop:cauchy_completion_is_cl}.
		\begin{proof} Consider the following diagram in $\lip$.
			\[\begin{tikzcd}[row sep = 1.5em, column sep = 1.5em]
				M\times M \arrow[rr, "\phi\times\phi"] \arrow[rd, "i\times i"]\arrow[ddd,"+_M"]& &N\times N\arrow[ddd,"+_N"]\\
				&  \overline{M}\times\overline{M}\arrow[ur,"\overline{\phi}\times\overline{\phi}"]\arrow[d,"+_{\overline{M}}"]\\
				&\overline{M}\arrow[dr,"\overline{\phi}"]&\\
				M\arrow[ur,"i"]\arrow[rr,"\phi"]&&N
			\end{tikzcd}\]
			The outer square commutes since $\phi$ is a monoid homomorphism, the upper and lower triangles commute by universality of the Cauchy completion, and the left trapezoid commutes by construction (Proposition \ref{prop:cauchy_completion_is_cl}). It follows that $+_N\circ (\overline{\phi}\times \overline{\phi})\circ (i\times i) = \overline{\phi}\circ +_{\overline{M}} \circ (i\times i)$. Since $i\times i$ is an epimorphism, it follows that $+_N\circ (\overline{\phi}\times \overline{\phi})= \overline{\phi}\circ +_{\overline{M}} $, which is precisely the statement that $\overline{\phi}$ respects the monoid operations. Moreover, $\overline{\phi}(0) = \overline{\phi}(i(0)) = \phi(0) = 0$ and thus $\overline{\phi}$ is a monoid homomorphism.
		\end{proof}
	\end{proposition}
	
	The following lemma shows that formal sums of Cauchy sequences are Cauchy sequences in the space of persistence diagrams.
	\begin{lemma}\label{lem:cauch_seq_in_diagrams}
		Let $(X,d,A)$ be a metric pair and let $(x_n^i)_n$ be Cauchy sequences in $X$ for $1\leq i\leq k$. Then the sequence $(x_n^1 + \dots + x_n^k)_n$ is a Cauchy sequence in $D(X,A)$.
		\begin{proof}
			Since the inclusion $(X,d)\hookrightarrow (D(X,A),W_p[d])$ is Lipschitz, each $(x_n^i)_n$ is a Cauchy sequence in $D(X,A)$, and since addition in $D(X,A)$ is Lipschitz, the term-wise sum of these sequences is also a Cauchy sequence in $D(X,A)$.
		\end{proof}
	\end{lemma}
	
	We now consider an extension of the space of finite persistence diagram considered earlier to a space of persistence diagrams allowing possibly countably infinitely many summands.
	\begin{definition}
		For a set $X$, let $\overline{D}(X) = \{f:X\to \N \ | \ f(x) = 0 \textup{ for all but countably many }x\in X\}$.
	\end{definition}
	Clearly $\overline{D}(X)$ is a commutative monoid with monoid operation given by pointwise addition of functions and identity the zero function. Moreover, if $A\subset X$ then we can identify $\overline{D}(A)$ with a submonoid of $\overline{D}(X)$. Elements of $\overline{D}(X)$ can also be viewed as countable formal sums of elements of $X$. There is a canonical map $i:X\to\overline{D}(X)$ which, when elements of $\overline{D}(X)$ are viewed as formal sums, is given by $x\mapsto x$. 
	
	\begin{definition} \label{def:overline-d}
		Let $(X,d,A)$ be a metric pair. Define the space of \textit{countable persistence diagrams on $X$ relative to $A$ by $\overline{D}(X,A)= \overline{D}(X)/\overline{D}(A)$}. We will also refer to elements of $\overline{D}(X,A)$ as \emph{countable persistence diagrams on $(X,A)$}.
	\end{definition}
	As with finite persistence diagrams, we have a monoid isomorphism $\overline{D}(X,A)\cong \overline{D}(X \setminus A)$. Effectively, the existence of this isomorphism shows that every element of the quotient $\overline{D}(X)/\overline{D}(A)$ has a unique representative which takes value zero on all of $A$, i.e., a formal sum with values in $X \setminus A$. We identify $D(X,A)$ as a submonoid of $\overline{D}(X,A)$ in the obvious way.
	
	The Wasserstein distances are easily extended to $\overline{D}(X,A)$. 
	\begin{definition}
		Let $\alpha,\beta\in \overline{D}(X,A)$. We call $\sigma \in \overline{D}(X\times X)$ a \textit{matching between $\alpha$ and $\beta$} if $\sum_{z\in X}\sigma(x,z) =\alpha(x)$ and $\sum_{z\in X} \sigma(z,x) = \beta(x)$ for all $x\in X\backslash A$. For $p\in [1,\infty]$, the \textit{$p$-cost} of $\sigma$ is defined by $\elp{(d(x_i,x_i'))_{i\in I}}{p}$, where $I$ is a countable set such that $\sigma = \sum_{i\in I} (x_i,x_i')$. The \textit{$p$-Wasserstein distance between $\alpha$ and $\beta$} is then defined to be 
		\[W_p[d](\alpha,\beta) = \inf_{\sigma}\textup{Cost}_p[d](\sigma),\]
		the infimum being taken over all matchings between $\alpha$ and $\beta$.
	\end{definition}
	
	Like the finite Wasserstein distance, $W_p$ is $p$-subadditive and translation sub-invariant, the proofs being essentially identical to those given in the finite case (see \cite[Lemma 4.17]{bubenik2019universality} and Lemma \ref{lem:translation_subinvariance} for the respective proofs in the finite case).
	
	In order to understand the Cauchy completion of the space of finite persistence diagrams in our framework, we will consider the following subspace of $\overline{D}(X,A)$. First, we need a preliminary definition. For $\alpha\in \overline{D}(X,A)$, let $\alpha-D(X,A) = \{\beta\in \overline{D}(X,A) \ | \ \alpha = \beta + \gamma \textup{ for some } \gamma\in D(X,A)\}$. In words, $\alpha-D(X,A)$ is the set of all countable persistence diagrams on $(X,A)$ which differ from $\alpha$ by a finite persistence diagram on $(X,A)$.
	
	\begin{definition}\label{def:diagrams_with_finite_distance_to_A}
		Let $(X,d,A)$ be a metric pair and let $p\in [1,\infty)$. We define 
		\begin{equation}\label{eq:p_finite_diagrams} \overline{D}_p(X,A) = \{ \alpha\in \overline{D}(X,A) \ | \ W_p[d](\beta,0)<\infty \textup{ for some } \beta\in \alpha - D(X,A)\}.
		\end{equation}
		For each $\epsilon > 0$ 
		let $U_{\epsilon} = \{x\in X \ | \ d(x,A) >\epsilon\}$. For $p = \infty$ we define
		\begin{equation}\label{eq:infinity_finite_diagrams}\overline{D}_\infty(X,A) = \{\alpha\in \overline{D}(X,A) \ | \ \sum_{x\in U_{\epsilon}} \alpha(x) <\infty \textup{ for all } \epsilon >0\}.
		\end{equation}
	\end{definition}
	For $p\in [1,\infty)$, elements of $\overline{D}_p(X,A)$ should be thought of as those countable persistence diagrams on $(X,A)$ which, after removing at most finitely many summands, have finite Wasserstein distance to the zero diagram. Elements of $\overline{D}_\infty(X,A)$ are those diagrams for which there are only finitely many summands at a given positive distance away from $A$. Evidently, $\overline{D}_p(X,A)$ is a submonoid of $\overline{D}(X,A)$ for all $p\in[1,\infty]$. The inclusion of $D(X,A)$ into $\overline{D}(X,A)$ is in fact an inclusion $D(X,A) \subset \overline{D}_p(X,A)$. The following lemma shows that $\overline{D}_p(X,A)\subset \overline{D}_\infty(X,A)$ for all $p\geq 1$.

\begin{lemma}\label{lem:p_finite_implies_infinity_finite}
		Let $(X,d,A)$ be a metric pair and let $p\in[1,\infty]$. Then $\overline{D}_p(X,A)\subset\overline{D}_\infty(X,A)$.
\begin{proof}
  Let $p\in[1,\infty)$ and let $\alpha\in \overline{D}_p(X,A)$.
  Then there is some $\beta\in \alpha-D(X,A)$ such that $W_p(\beta,0)<\infty$.
  Let $\eps > 0$.
Since $W_p(\beta,0) \geq ((\sum_{x \in U_{\epsilon}}\beta(x))\eps^p)^{\frac{1}{p}}$,
   $\sum_{x\in U_{\epsilon}}\beta(x)<\infty$. Since $\beta\in \alpha-D(X,A)$,  there is some $\gamma\in D(X,A)$ such that $\alpha = \beta + \gamma$. Since $\sum_{x \in X \setminus A} \gamma(x) < \infty$, $\sum_{x\in U_{\epsilon}}\gamma(x)<\infty$. Therefore $\sum_{x\in U_{\epsilon}}\alpha(x) = \sum_{x \in U_{\epsilon}} \beta(x) + \sum_{x \in U_{\epsilon}} \gamma(x) < \infty$ and thus $\alpha \in \overline{D}_\infty(X,A)$.
\end{proof} 
\end{lemma}

The following result was first proved for classical persistence diagrams for $p \in [1,\infty)$~\cite{mileyko2011probability} and for bar codes for $p=\infty$~\cite{blumberg2014robust}. With a few changes, these proofs generalize to metric pairs~\cite{irynapeter}.
	
	\begin{proposition}[{\cite{irynapeter}}]\label{prop:diags_complete_iff_X_complete}
		Let $(X,d,A)$ be a metric pair and fix $p\in [1,\infty]$.  Then $(\overline{D}_p(X,A),W_p[d])$ is complete if and only if the quotient space $(X/A,\overline{d}_1)$ is complete.
		%
		%
	\end{proposition}
This result will allow us to completely describe the Cauchy completion of the space of persistence diagrams on a metric pair $(X,d,A)$.
We will need the following lemmas.
\begin{lemma}\label{lem:infinite_wasserstein_as_limit}
  Let $(X,d,A)$ be a metric pair, with Cauchy completion $(\overline{X},\overline{d},A/\!\!\sim_0)$. Let $p\in [1,\infty]$. Consider the spaces $(D(X,A),W_p[d])$ and $(\overline{D}_p(\overline{X},A/\!\!\sim_0),W_p[\overline{d}])$. Then for ${\textstyle\sum_{j = 1}^\infty [x_n^j]},{\textstyle\sum_{i = j}^\infty [y_n^j]}\in \overline{D}_p(\overline{X},A/\!\!\sim_0)$, we have
		\[W_p[\overline{d}]({\textstyle\sum_{j = 1}^\infty [x_n^j]},{\textstyle\sum_{j = 1}^\infty [y_n^j]}) = \lim_{N\to \infty}\lim_{m\to\infty}W_p[d]({\textstyle\sum_{j = 1}^Nx_m^j},{\textstyle\sum_{j = 1}^Ny_m^j}).\]
\begin{proof}
Let ${\textstyle\sum_{j = 1}^\infty [x_n^j]},{\textstyle\sum_{j = 1}^\infty [y_n^j]}\in \overline{D}_p(\overline{X},A/\!\!\sim_0)$ be given. By translation sub-invariance, we have \[W_p[\overline{d}]({\textstyle\sum_{j = 1}^\infty [x_n^j]},{\textstyle\sum_{j = 1}^N [x_n^j]})\leq W_p[\overline{d}]({\textstyle\sum_{j = N+1}^\infty [x_n^j]},0).\]
By the definition of $\overline{D}_p(\overline{X},A/\!\!\sim_0)$, there is some $\beta \in {\textstyle\sum_{j = 1}^\infty [x_n^j]} - D(\overline{X},A/\!\!\sim_0)$ such that $W_p[\overline{d}](\beta,0) <\infty$. Hence we have that
$W_p[\overline{d}](\sum_{j=M}^{\infty}[x_n^j],0) < \infty$ for sufficiently large $M$. Therefore, for each $\epsilon >0$,
$W_p[\overline{d}]({\textstyle\sum_{j = N+1}^\infty [x_n^j]},0) <\epsilon$ for sufficiently large $N$ and hence $W_p[\overline{d}]({\textstyle\sum_{j = N+1}^\infty [x_n^j]},0)\to 0$ as $N\to\infty$. Therefore $\sum_{i = j}^N [x_n^j]\to\sum_{j = 1}^\infty [x_n^j]$ as $N\to\infty$.
			
			Now note that the canonical map $i:X\to \overline{X}$ induces a distance preserving map $(D(X,A),W_p[d])\to (\overline{D}(\overline{X},A/\!\!\sim_0),W_p[\overline{d}])$. Moreover, $\lim_{m\to\infty}i(x_m^j)=[x_n^j]$ for all $j\in \N$ (Lemma \ref{lem:limit_in_cauchy_completion}). Then
			\begin{multline*}
				 \lim_{N\to \infty}\lim_{m\to\infty}W_p[d]({\textstyle\sum_{j = 1}^Nx_m^j},{\textstyle\sum_{j = 1}^Ny_m^j})  =  \lim_{N\to \infty}\lim_{m\to\infty}W_p[\overline{d}]({\textstyle\sum_{j = 1}^Ni(x_m^j)},{\textstyle\sum_{i = 1}^Ni(y_m^j)})\\
				= W_p[\overline{d}]( \lim_{N\to \infty}\lim_{m\to\infty}{\textstyle\sum_{j = 1}^Ni(x_m^j)}, \lim_{N\to \infty}\lim_{m\to\infty}{\textstyle\sum_{j = 1}^Ni(y_m^j)})
				= W_p[\overline{d}]({\textstyle\sum_{j = 1}^\infty [x_n^j]},{\textstyle\sum_{j = 1}^\infty [y_n^j]}),
			\end{multline*}
			the second equality coming from the continuity of the metric. 
                        This completes the proof.
		\end{proof}
	\end{lemma}

For $\alpha\in \overline{D}_p(X,A)$ and $\epsilon>0$, let $U_\epsilon(\alpha)$ be given by $U_\epsilon(\alpha)(x) = \alpha(x)$ if $d(x,A)>\epsilon$ and $U_\epsilon(\alpha)(x) = 0$ otherwise. That is, $U_\epsilon(\alpha) = \sum_{x\in U_\epsilon}\alpha(x)$. By Lemma \ref{lem:p_finite_implies_infinity_finite}, $U_\epsilon(\alpha)\in D(X,A)$. Let $L_\epsilon(\alpha)\in\overline{D}_p(X,A)$ be given by $L_\epsilon(\alpha)(x) = \alpha(x)$ if $d(x,A)\leq \epsilon$ and $L_\epsilon(\alpha)(x) = 0$ otherwise.
	
	Note that $\alpha = L_\epsilon(\alpha) + U_\epsilon(\alpha)$. Moreover, by definition of $\overline{D}_p(X,A)$, we have that $W_p(L_\epsilon(\alpha),0)\to 0$ as $\epsilon\to 0$. Thus, by translation subinvariance, $W_p(\alpha,U_\epsilon(\alpha)) = W_p(L_\epsilon(\alpha) + U_\epsilon(\alpha),U_\epsilon(\alpha)) \leq W_p(L_\epsilon(\alpha),0)$ and hence $W_p(\alpha,U_\epsilon(\alpha))\to 0$ as $\epsilon\to 0$. Moreover, it is straightforward to see that for $\alpha,\beta\in \overline{D}_p(X,A)$, $\alpha = \beta$ if and only if $U_\epsilon(\alpha) = U_\epsilon(\beta)$ for all $\epsilon>0$.

	\begin{lemma}\label{lem:diagrams_on_sep_space_satisfies_separation}
		Let $(X,d,A)$ be a separable metric pair and let $p\in[1,\infty]$. The CL monoid $(\overline{D}_p(X,A),W_p[d])$ satisfies the separation axiom.
\begin{proof}
Let $\alpha,\beta\in\overline{D}_p(X,A)$ and suppose that $W_p(\alpha,\beta) = 0$. Then 
\begin{multline*}
  W_p(U_\epsilon(\alpha),U_\epsilon(\beta)) 
  \leq W_p(U_\epsilon(\alpha),\alpha)  + W_p(\alpha,\beta)  + W_p(\beta,U_\epsilon(\beta))\\ = W_p(U_\epsilon(\alpha),\alpha)  +  W_p(\beta,U_\epsilon(\beta))
\end{multline*} and hence $W_p(U_\epsilon(\alpha),U_\epsilon(\beta)) \to 0$ as $\epsilon\to 0$. Suppose that $U_{\epsilon_0}(\alpha) \neq U_{\epsilon_0}(\beta)$ for some $\epsilon_0>0$. Then there is exists an $x\in X$ with $d(x,A)>\epsilon_0$ such that $\alpha(x)\neq \beta(x)$. Suppose without loss of generality that $\alpha(x)>\beta(x)$.
For simplicity, assume $\alpha(x)=1$ and $\beta(x)=0$. The general argument is similar.
Let $\delta = \min(\min\{d(x,y) \ | \ y\neq x,\, \beta(y) \neq 0\}, d(x,A)-\epsilon_0)$. We claim that $W_\infty(\alpha,\beta)\geq \delta$. To see this, suppose that $W_\infty(\alpha,\beta)<\delta$, and let $\sigma\in \overline{D}(X\times X)$ be a matching between $\alpha$ and $\beta$ with $\textup{Cost}_\infty(\sigma)<\delta$. Then, since $d(x,A) \geq \delta+\epsilon_0 > \delta$,
there exists
$y\in X$ with $\beta(y) > 0$ and
$\sigma(x,y) > 0$. Therefore
$d(x,y)<\delta$. But no such point can exist by definition of $\delta$, and we have reached a contradiction. Thus $W_p(\alpha,\beta)\geq W_\infty(\alpha,\beta)\geq \delta >0$. But now this contradicts our assumption that $W_p(\alpha,\beta) = 0$. Hence $U_\epsilon(\alpha) = U_\epsilon(\beta)$ for all $\epsilon>0$ and thus $\alpha = \beta$.
		\end{proof}
              \end{lemma}
          
\begin{corollary}\label{cor:cauchy_completion_of_diagrams_satisfies_separation}
	Let $(X,d,A)$ be a metric pair and let $p\in [1,\infty]$. Then the CL monoid $(\overline{D}_p(\overline{X},A/{\sim_0}), W_p[\overline{d}])$ satisfies the separation axiom.
	\begin{proof}
		By the definition of the Cauchy completion, $(\overline{X},\overline{d},A/{\sim_0})$ satisfies the separation axiom, and so the result follows from Lemma \ref{lem:diagrams_on_sep_space_satisfies_separation}.
	\end{proof}
\end{corollary}

	\begin{theorem} \label{thm:cauchy}
		Let $p\in [1,\infty]$. Then $(\overline{D}_p(\overline{X},A/\!\!\sim_0),W_p[\overline{d}])$ is the Cauchy completion of $(D(X,A),W_p[d])$.
		\begin{proof}
			We will use Corollary \ref{cor:cauchy_completion_unique_upto_isom} to establish the result. Let $i:X\to \overline{X}$ denote the canonical map and define $j:D(X,A)\to \overline{D}_p(\overline{X},A/\!\!\sim_0)$ by $\sum_{k = 1}^M x^k \mapsto \sum_{k = 1}^M i(x^k)$. 
			We will verify that $(\overline{D}_p(\overline{X},A/\!\!\sim_0),W_p[\overline{d}])$ together with $j$ satisfy the universal property of the Cauchy completion (Theorem \ref{thm:univ_prop_of_cauchy_completion}). Since $\overline{X}$ is complete, so is $\overline{X}/(A/\!\!\sim_0)$. Thus, by Proposition \ref{prop:diags_complete_iff_X_complete}, $(\overline{D}_p(\overline{X},A/\!\!\sim_0),W_p[\overline{d}])$ is complete. By Lemma \ref{cor:cauchy_completion_of_diagrams_satisfies_separation}, $(\overline{D}_p(\overline{X},A/\!\!\sim_0),W_p[\overline{d}])$ satisfies the separation axiom as well.
			
			Let $(N,\rho)$ be a complete metric space satisfying the separation axiom and let $\phi:D(X,A)\to N$ be Lipschitz. Define $\overline{\phi}:(\overline{D}_p(\overline{X},A/\!\!\sim_0),W_p[\overline{d}])\to N$ by setting
			\[ \overline{\phi}({\textstyle\sum_{k = 1}^\infty[x^k_n]}) = \lim_{M\to \infty}\lim_{n\to \infty}\phi({\textstyle\sum_{k = 1}^Mx_n^k}).\]
			To see that $\overline{\phi}$ is well-defined, first note that, by Lemma \ref{lem:cauch_seq_in_diagrams},  since each $(x_n^k)_n$ is a Cauchy sequence in $X$, for each $M$ we have that $(\sum_{k = 1}^M x_n^k)_n$ is a Cauchy sequence in $D(X,A)$. Since $\phi$ is Lipschitz, $(\phi({\textstyle\sum_{k = 1}^Mx_n^k}))_n$ is a Cauchy sequence in $N$ and hence converges by completeness of $N$. Further, suppose that, for each $k$, $[x_n^k] = [y_n^k]$. Then $\lim_{n\to \infty}d(x_n^k,y_n^k) = 0$ for all $k$, and for all $M$ we have,
			\begin{multline*}
				\rho(\lim_{n\to \infty}\phi({\textstyle\sum_{k = 1}^Mx_n^k}),\lim_{n\to \infty}\phi({\textstyle\sum_{k = 1}^My_n^k})) = \lim_{n\to \infty}\rho(\phi({\textstyle\sum_{k = 1}^Mx_n^k}),\phi({\textstyle\sum_{k = 1}^My_n^k}))\\
				\leq \lipnorm{\phi}\lim_{n\to \infty}W_p[d]({\textstyle\sum_{k = 1}^Mx_n^k},{\textstyle\sum_{k = 1}^My_n^k})\leq \lipnorm{\phi}\lim_{n\to\infty}\textstyle\sum_{k = 1}^Md(x_n^k,y_n^k) = 0.
			\end{multline*}
			Since $(N,\rho)$ satisfies the separation axiom, we thus have $\lim_{n\to \infty}\phi(\sum_{k = 1}^Mx_n^k) = \lim_{n\to \infty}\phi(\sum_{k = 1}^My_n^k)$. Finally, consider the sequence $\left(\lim_{n\to \infty}\phi(\sum_{k = 1}^Mx_n^k)\right)_M$ in $N$. For an ordered pair $M<K$ we have
			\begin{multline*}
				\rho(\lim_{n\to \infty}\phi({\textstyle\sum_{k = 1}^Mx_n^k}),\lim_{n\to \infty}\phi({\textstyle\sum_{k = 1}^K x_n^k})) = \lim_{n\to \infty}\rho(\phi({\textstyle\sum_{k = 1}^Mx_n^k}),\phi({\textstyle\sum_{k = 1}^K x_n^k}))\\
				\leq \lipnorm{\phi}\lim_{n\to \infty}W_p[d]({\textstyle\sum_{k = 1}^Mx_n^k},{\textstyle\sum_{k = 1}^K x_n^k})\leq \lipnorm{\phi}\lim_{n\to \infty}\elp{(d(x_n^k,A))_{k = M+1}^K}{p}.
			\end{multline*}
			Since $\sum_{k = 1}^\infty [x^k_n]\in \overline{D}_p(\overline{X},A/\!\!\sim_0)$,
                        $\norm{(d(x_n^k,A))_{k=L}^{\infty}}_p < \infty$ for sufficiently large $L$. Hence
                        
            \[\lim_{n\to\infty}\elp{(d(x_n^k,A))_{k = M+1}^K}{p} = \elp{(\lim_{n\to\infty}d(x_n^k,A))_{k = M+1}^K}{p} = \elp{(\overline{d}([x_n^k],A/{\sim_0}))_{k = M+1}^K}{p} \to 0,\] as $M,K\to \infty$.
            Thus $\left(\lim_{n\to \infty}\phi(\sum_{k = 1}^Mx_n^k)\right)_{M}$  is a  Cauchy sequence in $N$ and hence converges, which shows that $\overline{\phi}$ is well-defined.
			
			Now $\overline{\phi}(j(x_1 + \dots + x_n)) = \overline{\phi}(i(x_1) + \dots + i(x_n)) = \phi(x_1 + \dots + x_n)$ since each $i(x_k)$ is just the constant Cauchy sequence on $x_k$. Thus $\overline{\phi}\circ j = \phi$. To see that $\overline{\phi}$ is unique, suppose that $\psi: (D_p(\overline{X},A/\sim_0)\to N$ is a Lipschitz map such that $\psi\circ j = \phi$. Then \begin{multline*}
					\psi({\textstyle\sum_{k = 1}^\infty[x^k_n]}) = \psi(\lim_{M\to\infty}{\textstyle\sum_{k = 1}^M[x^k_n]}) = \lim_{M\to\infty}\psi({\textstyle\sum_{k = 1}^M[x^k_n]})\\
					= \lim_{M\to\infty}\psi(\lim_{n\to\infty}{\textstyle\sum_{k = 1}^Mi(x^k_n)}) =\lim_{M\to\infty}\lim_{n\to \infty}\psi({\textstyle\sum_{k = 1}^Mi(x^k_n)}) \\
					= \lim_{M\to\infty}\lim_{n\to\infty}\psi({\textstyle j(\sum_{k = 1}^Mx^k_n)}) = \lim_{M\to\infty}\lim_{n\to \infty}\phi({\textstyle\sum_{k = 1}^Mx_n^k}) = \overline{\phi}({\textstyle\sum_{i = 1}^\infty[x^i_n]}),
				\end{multline*}				
			and thus $\psi = \phi$.			
			
			Lastly, we will show that $\lipnorm{\overline{\phi}} = \lipnorm{\phi}$. We have
			\begin{multline*}
				\rho( \overline{\phi}({\textstyle\sum_{i = 1}^\infty[x^i]_n}), \overline{\phi}({\textstyle\sum_{i = 1}^\infty[y^i]_n}) = \rho(\lim_{N\to \infty}\lim_{n\to \infty}\phi({\textstyle\sum_{i = 1}^Nx_n^i}),\lim_{N\to \infty}\lim_{n\to \infty}\phi({\textstyle\sum_{i = 1}^Ny_n^i}))\\
				= \lim_{n\to \infty}\lim_{N\to \infty} \rho(\phi({\textstyle\sum_{i = 1}^Nx_n^i}),\phi({\textstyle{\sum_{i = 1}^Ny_n^i}})) \leq \lipnorm{\phi} \lim_{n\to \infty}\lim_{N\to \infty} W_p[d]({\textstyle\sum_{i = 1}^Nx_n^i},\textstyle{\sum_{i = 1}^Ny_n^i})\\
				= \lipnorm{\phi}W_p[\overline{d}]({\textstyle\sum_{i = 1}^\infty[x^i]_n},{\textstyle\sum_{i = 1}^\infty[y^i]_n}),
			\end{multline*}
			the last equality coming from Lemma \ref{lem:infinite_wasserstein_as_limit}. Thus $\lipnorm{\overline{\phi}}\leq \lipnorm{\phi}$. On the other hand, we have $\lipnorm{\phi} = \lipnorm{\overline{\phi}\circ j}\leq \lipnorm{\overline{\phi}}\lipnorm{j}= \lipnorm{\overline{\phi}}$ since $\lipnorm{j} = 1$. Thus $\lipnorm{\overline{\phi}} = \lipnorm{\phi}$, completing the proof.		
		\end{proof}
	\end{theorem}

	\subsection{Universal property for the Cauchy completion of diagrams}
	In this section we extend the universality result of \cite{bubenik2019universality} to the Cauchy completion of the space of persistence diagrams.
	
	\begin{theorem}[Universal Property for {$(\overline{D}_p(\overline{X},A),W_p[\overline{d}],+)$}] \label{thm:univ_prop_for_cauchy_completion_of_diagrams} Let $(X,d,A)$ be a metric pair and let $p\in [1,\infty]$. Let $(N,\rho, +, 0)$ be a complete CL monoid satisfying the separation axiom, and let $\phi:X\to N$ be Lipschitz and satisfy $\phi(A) = \{0\}$. Then 
		\begin{thmenum}
			\item \label{thm:univ_prop_for_cauchy_completion_of_diagram_part_1} There is a unique Lipschitz monoid homomorphism $\tilde{\phi}:(\overline{D}_p(\overline{X},A), W_p[\overline{d}])\to (N,\rho)$ such that $\phi = \tilde{\phi}\circ \iota$, where $\iota:X\to \overline{D}_p(\overline{X},A)$ denotes the composition $X\to D(X,A)\to \overline{D}_p(\overline{X},A)$; and
			\item \label{thm:univ_prop_for_cauchy_completion_of_diagram_part_2}
			$\lipnorm{\tilde{\phi}} = \lipnorm{\phi}$.
		\end{thmenum}
		\begin{proof}
			Let $\delta:X\to D(X,A)$ denote the canonical inclusion of $X$ to its space of persistence diagrams and let $i:D(X,A)\to \overline{D}_p(\overline{X},A)$ denote the canonical map of $D(X,A)$ to its Cauchy completion. By the universal property of the space of persistence diagrams~\cite{bubenik2019universality}, there is a unique monoid homomorphism $\phi':D(X,A)\to N$ with $\phi'(A) = \{0\}$ and such that $\phi'\circ \delta = \phi$. Then by the universal property of the Cauchy completion (Theorem \ref{thm:univ_prop_of_cauchy_completion}), there is a unique map Lipschitz map  $\tilde{\phi}:\overline{D}_p(\overline{X},A)\to N$ such that $\tilde{\phi}\circ i = \phi'$. Then $\tilde{\phi}\circ \iota = \tilde{\phi}\circ i\circ \delta = \phi'\circ \delta = \phi$ and $\lipnorm{\tilde{\phi}} = \lipnorm{\phi'} = \lipnorm{\phi}$. Moreover, by Proposition \ref{prop:unique_morphism_is_monoid_homomorphism}, $\tilde{\phi}$ is a monoid homomorphism, completing the proof.
		\end{proof}	
	\end{theorem}

\section{Universal Banach spaces} \label{sec:universal-banach}

We will show that every pointed metric space has a canonical embedding into a Banach space.
This universal Banach space has been constructed independently many times \cite{arens1956embedding, flood1984free, pestov1986free, godefroy2003lipschitz-free,weaver2018lipschitz}. It is sometimes referred to as the \emph{Arens-Eells space} \cite{weaver2018lipschitz} or the Lipschitz-free Banach space \cite{godefroy2003lipschitz-free}. We present a construction based on the partial optimal transport problem, which we formulate in the language of linear programming.

In this section we restrict to metric spaces $(X,d)$ that satisfy the separation and finiteness conditions. That is, our notion of metric spaces agrees with the usual one. In particular, $\cat{Lip_*}$ will denote the category whose objects are pointed metric spaces satisfying the separation and finiteness conditions and whose morphisms are pointed Lipschitz maps.
Also, all vector spaces will be real vector spaces.

\subsection{Normed vector spaces} \label{sec:nvs}

Let $\cat{NVS}$ denote the category whose objects are normed vector spaces $(V,\norm{\ })$ and whose morphisms are bounded linear operators, $T:(V,\norm{ }_V) \to (W,\norm{\ }_W)$, i.e. linear maps $T:V \to W$ such that there exists a $K$ for which $\norm{Tv}_W \leq K \norm{v}_V$ for all $v \in V$.
For a bounded linear operator $T:V \to W$, let $\norm{T}_{\op}$ denote the operator norm of $T$, which equals the infimum of all $K$ such that $\norm{Tv}_W \leq K\norm{v}_V$ for all $v \in V$.
There is a forgetful functor $U: \cat{NVS} \to \cat{Lip_*}$ given by $U(V,\norm{\ }) = (V,d_{\norm{\ }},0)$, where $d_{\norm{\ }}(v,w) = \norm{v-w}$, and
$U(T:(V,\norm{\ }_V) \to (W,\norm{\ }_W)) = T:(V,d_{\norm{\ }_V},0) \to (W,d_{\norm{\ }_W},0)$.
Note that $\norm{T}_{\Lip} = \norm{T}_{\op}$.

Given a set $X$, let $V(X)$ denote the free vector space on $X$ of formal $\R$-linear combinations of elements of $X$.
That is,
$V(X)$ is the set of functions from $X$ to $\R$ with finite support, together with pointwise addition and scalar multiplication.
By identifying $x \in X$ with the indicator function on $x$, we may write each $\mu \in V(X)$ as $\mu = \sum_{i=1}^n \mu_i x_i$ for some $n \geq 0$, $\mu_i \in \R \setminus 0$, and distinct $x_i \in X$.
For a pointed set $(X,x_0)$, let $V(X,x_0) = V(X)/V(x_0)$.
Define $\rho:V(X) \to V(X \setminus x_0)$ by
$\rho(\mu) = \mu|_{X \setminus x_0}$, or equivalently, 
$\rho(\sum_{i=0}^n \mu_i x_i) = \sum_{i=1}^n \mu_i x_i$.
Then $\rho$ induces a well-defined map $\overline{\rho}: V(X,x_0) \to V(X \setminus x_0)$.
Define $\iota: V(X \setminus x_0) \to V(X)$ by $\iota(\sum_{i=1}^n \mu_i x_i) = \sum_{i=0}^n \mu_i x_i$, where $\mu_0 = -\sum_{i=1}^n \mu_i$.
Then $\iota$ induces a map $\overline{\iota}:V(X \setminus x_0) \to V(X,x_0)$ and $\overline{\rho}$ and $\overline{\iota}$ are inverse maps.
Define $r: V(X,x_0) \to V(X)$ by $r = \iota \circ \overline{\rho}$.
Then
for $\mu \in V(X,x_0)$,
$r$ chooses a canonical representative $\sum_{i=0}^n \mu_i x_i$ with $\sum_{i=0}^n \mu_i = 0$.
A set map $f:X \to Y$ extends to a linear map $f:V(X) \to V(Y)$ by $f(\sum_{i=1}^n \mu_i x_i) = \sum_{i=1}^n \mu_i f(x_i)$.
Furthermore, a pointed set map $f:(X,x_0) \to (Y,y_0)$ extends to a linear map $f:V(X,x_0) \to V(Y,y_0)$.

\begin{definition} \label{def:wasserstein-norm}
  Given $(X,d,x_0) \in \cat{Lip_*}$, we define the \emph{Wasserstein norm} on $V(X,x_0)$ as follows. 
  For $\mu \in V(X,x_0)$, consider $r(\mu) = \sum_{i=0}^n \mu_i x_i$ with $\sum_{i=0}^n \mu_i = 0$. Define $\norm{\mu}_{W_1[d]}$ to be the solution of the linear programming problem
  \begin{align} \label{eq:wasserstein-norm}
    \begin{split}
    \text{minimize }
    & \sum_{i=0}^n \sum_{j=0}^n \pi_{ij} d(x_i,x_j)\\
    \text{subject to } & \sum_{k=0}^n (\pi_{ik} - \pi_{ki}) = \mu_i, i=0,\ldots,n\\
                               & \pi_{ij} \geq 0, i,j=0,\ldots,n.
    \end{split}
  \end{align}
  Call $\pi = \sum_{i=0}^n \sum_{j=0}^n \pi_{i,j}(x_i,x_j) \in V(X \times X)$ that satisfies the constraint of \eqref{eq:wasserstein-norm} a \emph{coupling} for $\mu$ and call $\sum_{i=0}^n \sum_{j=0}^n \pi_{ij} d(x_i,x_j)$ the \emph{cost} of the coupling, denote $\Cost_1[d](\pi)$.
\end{definition}

\begin{proposition} \label{prop:wasserstein-norm}
  For $(X,d,x_0) \in \cat{Lip_*}$, the Wasserstein norm $\norm{\ }_{W_1[d]}$ is a norm on $V(X,x_0)$.
\end{proposition}

\begin{proof}
  For subadditivity, consider $\mu, \mu' \in V(X,x_0)$ with $r(\mu) = \sum_{i=0}^n \mu_i x_i$ and $r(\mu') = \sum_{i=0}^{n'} \mu'_i x'_i$, where $x'_0 = x_0$.
  Let $\nu = \mu + \mu'$.
  Then $r(\nu) = \sum_{i=0}^{n+n'} \nu_i x_i$, where for $i=1,\ldots,n$, $\nu_i = \mu_i$, for $i=1,\ldots,n'$, $\nu_{n+i} = \mu'_i$ and $x_{n+i} = x'_i$, and $\nu_0 = \mu_0 + \mu'_0$.
  Let $(\pi_{ij})_{i,j=0}^n$ be a coupling for a solution to \eqref{eq:wasserstein-norm} for $\mu$.
  Let $(\pi'_{ij})_{i,j=0}^n$ be a coupling for a solution to \eqref{eq:wasserstein-norm} for $\mu'$.
  Define
  \[
  \tau = 
  \begin{bmatrix}
    \pi_{00}+\pi'_{00} & \pi_{01} & \cdots & \pi_{0n} & \pi'_{01} & \cdots & \pi'_{0n'}\\
    \pi_{10} & \pi_{11} & \cdots & \pi_{1n} & 0 & \cdots & 0 \\
    \vdots & \vdots & \ddots & \vdots & \vdots & \ddots & \vdots \\
    \pi_{n0} & \pi_{n1} & \cdots & \pi_{nn} & 0 & \cdots & 0 \\
    \pi'_{10} & 0 & \cdots & 0 & \pi'_{11} & \cdots & \pi'_{1n'}\\
    \vdots & \vdots & \ddots & \vdots & \vdots & \ddots & \vdots \\
    \pi'_{n'0} & 0 & \cdots & 0 & \pi'_{n'1} & \cdots & \pi'_{n'n'}\\
\end{bmatrix}.
\]
Then $\tau$ is a coupling for $\nu$ and its cost is the sum of the costs of $\pi$ and $\pi'$.
Therefore $\norm{\mu + \mu'}_{W_1[d]} \leq \norm{\mu}_{W_1[d]} + \norm{\mu'}_{W_1[d]}$.

Let $\mu \in V(X,x_0)$ and $a > 0$. Then $\pi$ is a coupling for $\mu$ if and only if $a\pi$ is a coupling for $a\mu$. Furthermore $\pi$ is a coupling for $\mu$ if and only if $\pi^\top$, the transpose of $\pi$, is a coupling for $-\mu$.
Since $\Cost_1[d](a\pi) = a\Cost_1[d](\pi)$ and $\Cost_1[d](\pi^{\top}) = \Cost_1[d](\pi)$, it follows that for all $b \in \R$, $\norm{b\mu}_{W_1[d]} = \abs{b}\norm{\mu}_{W_1[d]}$.

Finally, let $\mu \in V(X,x_0)$ with $r(\mu) = \sum_{i=0}^n \mu_i x_i$ and $\norm{\mu}_{W_1[d]} = 0$.
Then there is a coupling $\pi$ for $\mu$ with zero cost.
Thus, by the objective function in \eqref{eq:wasserstein-norm}, $\pi_{ij}=0$ for all $i\neq j$.
Therefore, by the constraint in \eqref{eq:wasserstein-norm}, $\mu_i = 0$ for all $i$ and hence $\mu = 0$.
\end{proof}

\begin{theorem} \label{thm:nvs-embedding}
  Let $(X,d,x_0) \in \cat{Lip_*}$.
  The canonical map $i: (X,d,x_0) \to (V(X,x_0),d_{\norm{\ }_{W_1[d]}}, 0)$ is an isometric embedding.
\end{theorem}

\begin{proof}
  We need to verify that for all $x_1,x_2 \in X$,
if $\mu$ is the equivalence class of $x_1-x_2$ in $V(X,x_0)$, then
  $\norm{\mu}_{W_1[d]} = d(x_1,x_2)$.
  To start, $r(\mu) = 0x_0 + 1x_1 + (-1)x_2$.
  Let $\pi$ be a coupling for $\mu$.
  Since the terms $\pi_{ii}$ have no impact on \eqref{eq:wasserstein-norm}, we will assume they equal $0$.
  The constraints in \eqref{eq:wasserstein-norm} give us the following.
  \begin{align}
    \pi_{01} + \pi_{02} - \pi_{10} - \pi_{20} &= 0 \label{eq:a1} \\
    \pi_{10} + \pi_{12} - \pi_{01} - \pi_{21} &= 1 \label{eq:a2} \\
    \pi_{20} + \pi_{21} - \pi_{02} - \pi_{12} &= -1 \nonumber 
  \end{align}
  Rearranging \eqref{eq:a1}, we obtain
  \begin{equation} \label{eq:a4}
    \pi_{10} = \pi_{01} + \pi_{02} - \pi_{20}.
  \end{equation}
  Substituting \eqref{eq:a4} in \eqref{eq:a2}, we obtain
  \begin{equation} \label{eq:a5}
    \pi_{12} = 1 + \pi_{21} + \pi_{20} - \pi_{02}.
  \end{equation}
  The cost of $\pi$ is given by
  \begin{equation*}
    (\pi_{12}+\pi_{21})d(x_1,x_2) + (\pi_{01}+\pi_{10})d(x_1,x_0) + (\pi_{02}+\pi_{20})d(x_2,x_0). 
  \end{equation*}
  Substituting \eqref{eq:a5} and \eqref{eq:a4}, we get
  \begin{equation*}
    (1+2\pi_{21}+\pi_{20}-\pi_{02})d(x_1,x_2) + (2\pi_{01}+\pi_{02}-\pi_{20})d(x_1,x_0) + (\pi_{02}+\pi_{20})d(x_2,x_0).
  \end{equation*}
  Rearranging we obtain
  \begin{multline} \label{eq:a6}
    d(x_1,x_2) + 2\pi_{21}d(x_1,x_2) + \pi_{20}[d(x_1,x_2) - d(x_1,x_0) + d(x_2,x_0)]\\ + 2\pi_{01}d(x_1,x_0) + \pi_{02}[d(x_1,x_0) + d(x_2,x_0) - d(x_1,x_2)].
  \end{multline}
  By the triangle inequality, the two sums inside the square brackets are both nonnegative.
  It follows that the cost is minimized by setting
  $\pi_{21} = \pi_{20} = \pi_{01} = \pi_{02} = 0$.
  By \eqref{eq:a4} and \eqref{eq:a5}, $\pi_{10} = 0$ and $\pi_{12}=1$, respectively.
  Furthermore, by \eqref{eq:a6}, the cost of this optimal coupling is $d(x_1,x_2)$.
\end{proof}

\begin{theorem} \label{thm:universal-nvs}
  Let $(X,d,x_0) \in \cat{Lip_*}$.
  For any $(W,\norm{\ }) \in \cat{NVS}$ and $\varphi: (X,d,x_0) \to (W,d_{\norm{\ }},0) \in \cat{Lip_*}$, there exists a unique map $\tilde{\varphi}: (V(X,x_0),\norm{\ }_{W_1[d]}) \to (W,\norm{\ }) \in \cat{NVS}$ such that $U\tilde{\varphi} \circ i = \varphi$.
  Furthermore, $\norm{U\tilde{\varphi}}_{\Lip} = \norm{\varphi}_{\Lip}$.
\end{theorem}

\begin{proof}
  Let $\mu \in V(X,x_0)$ with $r(\mu) = \sum_{i=0}^n \mu_i x_i$.
  Since $\tilde{\varphi}$ is a linear map, we must have that $\tilde{\varphi}(\mu) = \sum_{i=1}^n \mu_i \varphi(x_i)$ (since $\varphi$ is a pointed map, $\varphi(x_0) = 0$).

  Let $\pi$ be a coupling for $\mu$ which is a solution to the linear program \eqref{eq:wasserstein-norm}.
  Then $r(\mu) = \sum_{i=0}^n \sum_{j=0}^n (\pi_{ij} - \pi_{ji}) x_i
  = \sum_{i=0}^n \sum_{j=0}^n \pi_{ij} (x_i - x_j)$.
  Therefore
  $\tilde{\varphi}(\mu)
  = \sum_{i=0}^n \sum_{j=0}^n \pi_{ij} \tilde{\varphi} (x_i - x_j)
  = \sum_{i=0}^n \sum_{j=0}^n \pi_{ij} (\varphi(x_i) - \varphi(x_j))$.
  By the subadditivity of the norm $\norm{\ }$,
  \begin{align*}
  \norm{\tilde{\varphi}(\mu)}
  &\leq \sum_{i=0}^n \sum_{j=0}^n \pi_{ij} \norm{\varphi(x_i) - \varphi(x_j)}\\
  &\leq \norm{\varphi}_{\Lip} \sum_{i=0}^n \sum_{j=0}^n \pi_{ij} d(x_i,x_j)\\
  &= \norm{\varphi}_{\Lip} \norm{\mu}.
  \end{align*}
  Therefore $\tilde{\varphi}$ is a bounded linear map.
  Furthermore $\norm{U\tilde{\varphi}}_{\Lip} = \norm{\tilde{\varphi}}_{\op} = \norm{\varphi}_{\Lip}$.
\end{proof}



  

\begin{corollary}
  The forgetful functor $U: \cat{NVS} \to \cat{Lip_*}$ has a left adjoint $V: \cat{Lip_*} \to \cat{NVS}$ which sends $(X,d,x_0)$ to $(V(X,x_0),\norm{\ }_{W_1[d]})$ and sends $f:(X,d,x_0) \to (Y,d',y_0)$ to the induced linear map $f:(V(X,x_0),\norm{\ }_{W_1[d]}) \to (V(Y,y_0),\norm{\ }_{W_1[d']})$.
  Furthermore, $\norm{f}_{\op} = \norm{f}_{\Lip}$.
\end{corollary}

\subsection{Banach spaces} \label{sec:banach}

Recall that a Banach space is a complete normed vector space.
Let $\cat{Ban}$ denote the full subcategory of $\cat{NVS}$ consisting of Banach spaces.
Every normed vector space $V$ isometrically embeds as a dense vector subspace into its Cauchy completion $\hat{V}$.
This defines a functor $C: \cat{NVS} \to \cat{Ban}$ that is left adjoint to the inclusion functor $\cat{Ban} \hookrightarrow \cat{NVS}$.

Combining this with Theorem~\ref{thm:nvs-embedding}, we have the following.

\begin{corollary}
  Let $(X,d,x_0) \in \cat{Lip_*}$.
  The canonical map $i: (X,d,x_0) \to (\hat{V}(X,x_0),d_{\norm{\ }_{W_1[d]}}, 0)$ is an isometric embedding.
\end{corollary}

Since adjoint functors compose we have the following.

\begin{corollary}
  The forgetful functor $U: \cat{Ban} \to \cat{Lip_*}$ has a left adjoint $V: \cat{Lip_*} \to \cat{Ban}$ which sends $(X,d,x_0)$ to $(\hat{V}(X,x_0),\norm{\ }_{W_1[d]})$ and sends $f:(X,d,x_0) \to (Y,d',y_0)$ to the induced linear map $f:(\hat{V}(X,x_0),\norm{\ }_{W_1[d]}) \to (\hat{V}(Y,y_0),\norm{\ }_{W_1[d']})$.
  Furthermore, $\norm{f}_{\op} = \norm{f}_{\Lip}$.
\end{corollary}

Restating this as a universal property we have the following.

\begin{theorem} \label{thm:universal-banach}
  Let $(X,d,x_0) \in \cat{Lip_*}$.
  For any $(W,\norm{\ }) \in \cat{Ban}$ and $\varphi: (X,d,x_0) \to (W,d_{\norm{\ }},0) \in \cat{Lip_*}$, there exists a unique map $\tilde{\varphi}: (\hat{V}(X,x_0),\norm{\ }_{W_1[d]}) \to (W,\norm{\ }) \in \cat{Ban}$ such that $U\tilde{\varphi} \circ i = \varphi$.
  Furthermore $\norm{U\tilde{\varphi}}_{\Lip} = \norm{\varphi}_{\Lip}$.
\end{theorem}

\subsection{Isometric embeddings} \label{sec:embedding}

We summarize the embeddings we have constructed as follows.

\begin{theorem} \label{thm:isometric-embedding}
  Let $(X,d,x_0)$ be a pointed metric space. We have the following sequence of isometric embeddings of pointed metric spaces, where we omit the basepoint $0$ from the notation,
  \begin{equation*} \label{eq:isometric-embedding}
    (X,d,x_0) \incl (D(X,x_0),W_1[d]) \incl (K(X,x_0),W_1[d]) \incl (V(X,x_0),d_{\norm{\ }_{W_1[d]}}) \incl (\hat{V}(X,x_0),d_{\norm{\ }_{W_1[d]}})
  \end{equation*}
  and $(\overline{D}_1(X,x_0),W_1[d]) \incl (\hat{V}(X,x_0),d_{\norm{\ }_{W_1[d]}})$.
\end{theorem}

\begin{proof}
  The first map is an isometric embedding by \cite[Lemma 4.18]{bubenik2019universality}.
  The second map is an isometric embedding by Corollary~\ref{cor:vpd-w1}.
  It is well known that the last map is an isometric embedding.
  We prove that the third map is an isometric embedding.

  Let $\alpha,\beta,\gamma,\delta \in D(X,x_0)$.
  Then
  \begin{equation*}
    d_{\norm{\ }_{W_1[d]}}(\alpha-\beta,\gamma-\delta) = \norm{\alpha-\beta-(\gamma-\delta)}_{W_1[d]} =
\norm{\alpha+\delta-(\beta+\gamma)}_{W_1[d]}.
\end{equation*}
From \eqref{eq:wasserstein-norm}, we see that the latter equals the cost of the solution to the transshipment problem (Definition~\ref{def:transshipment}) from $\alpha+\delta$ to $\beta+\gamma$.
Since $p=1$, this equals the cost of the solution to the transportation problem (Definition~\ref{def:transportation_problem}) from $\alpha+\delta$ to $\beta+\gamma$.
That is, $\norm{\alpha+\delta-(\beta+\gamma)}_{W_1[d]} = W_1[d](\alpha+\delta,\beta+\gamma) = W_1[d](\alpha-\beta,\gamma-\delta)$.  

Since $(D(X,x_0),W_1[d])$ isometrically embeds in $(V(X,x_0),d_{\norm{\ }_{W_1[d]}})$, it follows that the Cauchy completion of the former isometrically embeds in the Cauchy completion of the latter.
\end{proof}

As a result, we may refer to each of the metrics in 
Theorem~\ref{thm:isometric-embedding}
as $W_1[d]$.

Given a metric pair $(X,d,A)$, define $V(X,A) = V(X)/V(A)$. By \cite[Remark 4.14]{bubenik2019universality}, for a metric pair $(X,d,A$), we have and isometric isomorphism $(D(X,A),W_1[d],+,0)\cong (D(X/A,A),W_1[\overline{d}])$. Combining this with Theorem \ref{thm:isometric-embedding}, we obtain the following.

\begin{corollary} \label{cor:isometric-embedding}
  Let $(X,d,A)$ be a metric pair. We have the following sequence of isometric embeddings of pointed metric spaces, where we omit the basepoint $0$ from the notation,
  \begin{equation*} 
    (X/A,\overline{d},A) \incl (D(X,A),W_1[d]) \incl (K(X,A),W_1[d]) \incl (V(X,A),W_1[d]) \incl (\hat{V}(X,A),W_1[d])
  \end{equation*}
  and $(\overline{D}_1(X,A),W_1[d]) \incl (\hat{V}(X,A),d_{\norm{\ }_{W_1[d]}})$.
\end{corollary}

\begin{example}
  Consider the metric pair $(\R^2_{\leq},d,\Delta)$. We have isometric embeddings of the spaces of classical persistence diagrams $(D(\R^2_{\leq},\Delta),W_1[d])$ and
$(\overline{D}_1(\R^2_{\leq},\Delta),W_1[d])$
  into the Banach space $(\hat{V}(\R^2,\Delta),\norm{\ }_{W_1[d]})$.
\end{example}
        
\subsection*{Acknowledgments}

The authors would like to thank Parker Edwards, Iryna Hartsock, Facundo M\'{e}moli, Luis Scoccola, and Alex Wagner for useful comments and beneficial conversations.
This research was partially supported by the Southeast Center for Mathematics and Biology, an NSF-Simons Research Center for Mathematics of Complex Biological Systems, under National Science Foundation Grant No.\ DMS-1764406 and Simons Foundation Grant No.\ 594594. This material is based upon work supported by, or in part by, the Army Research Laboratory and the Army Research Office under contract/grant number W911NF-18-1-0307.

\vspace{1ex} 

On behalf of all authors, the corresponding author states that there is no conflict of interest.
	

\begin{thebibliography}{CGGGS21}
 	
 	\bibitem[AE56]{arens1956embedding}
 	Richard~F. Arens and James Eells, Jr.
 	\newblock On embedding uniform and topological spaces.
 	\newblock {\em Pacific J. Math.}, 6:397--403, 1956.
 	
 	\bibitem[AENY19]{asashiba2019approximation}
 	Hideto Asashiba, Emerson~G Escolar, Ken Nakashima, and Michio Yoshiwaki.
 	\newblock On approximation of $2 $ d persistence modules by
 	interval-decomposables.
 	\newblock {\em arXiv preprint arXiv:1911.01637}, 2019.
 	
 	\bibitem[BBE21]{betthauser2019graded}
 	Leo Betthauser, Peter Bubenik, and Parker~B Edwards.
 	\newblock Graded persistence diagrams and persistence landscapes.
 	\newblock {\em Discrete \& Computational Geometry}, pages 1--28, 2021.
 	
 	\bibitem[BBI01]{bbi:book}
 	Dmitri Burago, Yuri Burago, and Sergei Ivanov.
 	\newblock {\em A course in metric geometry}, volume~33 of {\em Graduate Studies
 		in Mathematics}.
 	\newblock American Mathematical Society, Providence, RI, 2001.
 	
 	\bibitem[BE21]{bubenik2019universality}
 	Peter Bubenik and Alex Elchesen.
 	\newblock Universality of persistence diagrams and the bottleneck and
 	{W}asserstein distances.
 	\newblock {\em Computational Geometry, accepted, (arXiv:1912.02563)}, 2021.
 	
 	\bibitem[BGMP14]{blumberg2014robust}
 	Andrew~J Blumberg, Itamar Gal, Michael~A Mandell, and Matthew Pancia.
 	\newblock Robust statistics, hypothesis testing, and confidence intervals for
 	persistent homology on metric measure spaces.
 	\newblock {\em Foundations of Computational Mathematics}, 14(4):745--789, 2014.
 	
 	\bibitem[BH21]{irynapeter}
 	Peter Bubenik and Iryna Hartsock.
 	\newblock Topological and metric properties of spaces of generalized
 	persistence diagrams.
 	\newblock {\em arXiv preprint}, 2021.
 	
 	\bibitem[BOO21]{botnan2021signed}
 	Magnus~Bakke Botnan, Steffen Oppermann, and Steve Oudot.
 	\newblock Signed barcodes for multi-parameter persistence via rank
 	decompositions and rank-exact resolutions.
 	\newblock {\em arXiv preprint arXiv:2107.06800}, 2021.
 	
 	\bibitem[BSS18]{bubenik2018wasserstein}
 	Peter Bubenik, Jonathan Scott, and Donald Stanley.
 	\newblock Wasserstein distance for generalized persistence modules and abelian
 	categories.
 	\newblock {\em arXiv preprint arXiv:1809.09654}, 2018.
 	
 	\bibitem[CGGGS21]{Che:2021}
 	Mauricio Che, Fernando Galaz-Garc{\'\i}a, Luis Guijarro, and Ingrid~Membrillo
 	Solis.
 	\newblock Metric geometry of spaces of persistence diagrams.
 	\newblock 09 2021.
 	
 	\bibitem[CSEH07]{cohen2007stability}
 	David Cohen-Steiner, Herbert Edelsbrunner, and John Harer.
 	\newblock Stability of persistence diagrams.
 	\newblock {\em Discrete Comput. Geom.}, 37(1):103--120, 2007.
 	
 	\bibitem[CSEHM10]{cohen2010lipschitz}
 	David Cohen-Steiner, Herbert Edelsbrunner, John Harer, and Yuriy Mileyko.
 	\newblock Lipschitz functions have {$L_p$}-stable persistence.
 	\newblock {\em Found. Comput. Math.}, 10(2):127--139, 2010.
 	
 	\bibitem[CZCG04]{Collins:2004}
 	Anne Collins, Afra Zomorodian, Gunnar Carlsson, and Leonidas~J. Guibas.
 	\newblock A barcode shape descriptor for curve point cloud data.
 	\newblock {\em Computers \& Graphics}, 28(6):881 -- 894, 2004.
 	
 	\bibitem[DL21]{divol2019understanding}
 	Vincent Divol and Th\'{e}o Lacombe.
 	\newblock Understanding the topology and the geometry of the space of
 	persistence diagrams via optimal partial transport.
 	\newblock {\em J. Appl. Comput. Topol.}, 5(1):1--53, 2021.
 	
 	\bibitem[Edw11]{edwards2011KRtheorem}
 	D.~A. Edwards.
 	\newblock On the {K}antorovich-{R}ubinstein theorem.
 	\newblock {\em Expo. Math.}, 29(4):387--398, 2011.
 	
 	\bibitem[ELZ00]{edelsbrunner2000topological}
 	Herbert Edelsbrunner, David Letscher, and Afra Zomorodian.
 	\newblock Topological persistence and simplification.
 	\newblock In {\em Foundations of Computer Science, 2000. Proceedings. 41st
 		Annual Symposium on}, pages 454--463. IEEE, 2000.
 	
 	\bibitem[EM19]{elchesen2019reflection}
 	Alexander Elchesen and Facundo M\'{e}moli.
 	\newblock The reflection distance between zigzag persistence modules.
 	\newblock {\em J. Appl. Comput. Topol.}, 3(3):185--219, 2019.
 	
 	\bibitem[FG10]{FIGALLI2010107}
 	Alessio Figalli and Nicola Gigli.
 	\newblock A new transportation distance between non-negative measures, with
 	applications to gradients flows with dirichlet boundary conditions.
 	\newblock {\em Journal de Mathématiques Pures et Appliquées}, 94(2):107--130,
 	2010.
 	
 	\bibitem[Flo84]{flood1984free}
 	Joe Flood.
 	\newblock Free topological vector spaces.
 	\newblock {\em Dissertationes Math. (Rozprawy Mat.)}, 221:95, 1984.
 	
 	\bibitem[GK03]{godefroy2003lipschitz-free}
 	G.~Godefroy and N.~J. Kalton.
 	\newblock Lipschitz-free {B}anach spaces.
 	\newblock volume 159, pages 121--141. 2003.
 	\newblock Dedicated to Professor Aleksander Pe\l czy\'{n}ski on the occasion of
 	his 70th birthday.
 	
 	\bibitem[GL21]{giusti2021signatures}
 	Chad Giusti and Darrick Lee.
 	\newblock Signatures, lipschitz-free spaces, and paths of persistence diagrams.
 	\newblock {\em arXiv preprint arXiv:2108.02727}, 2021.
 	
 	\bibitem[Hun80]{hungerford_algebra}
 	Thomas~W. Hungerford.
 	\newblock {\em Algebra}, volume~73 of {\em Graduate Texts in Mathematics}.
 	\newblock Springer-Verlag, New York-Berlin, 1980.
 	\newblock Reprint of the 1974 original.
 	
 	\bibitem[Kel85]{kellerer1982duality}
 	Hans~G. Kellerer.
 	\newblock Duality theorems and probability metrics.
 	\newblock In {\em Proceedings of the seventh conference on probability theory
 		({B}ra\c{s}ov, 1982)}, pages 211--220. VNU Sci. Press, Utrecht, 1985.
 	
 	\bibitem[KM21]{kim2021generalized}
 	Woojin Kim and Facundo M{\'e}moli.
 	\newblock Generalized persistence diagrams for persistence modules over posets.
 	\newblock {\em Journal of Applied and Computational Topology}, 5(4):533--581,
 	2021.
 	
 	\bibitem[Mai12]{mainini2011}
 	Edoardo Mainini.
 	\newblock A description of transport cost for signed measures.
 	\newblock {\em Journal of Mathematical Sciences}, 181(6):837--855, 2012.
 	
 	\bibitem[MMH11]{mileyko2011probability}
 	Yuriy Mileyko, Sayan Mukherjee, and John Harer.
 	\newblock Probability measures on the space of persistence diagrams.
 	\newblock {\em Inverse Problems}, 27(12):124007, 2011.
 	
 	\bibitem[MP20]{mccleary2020edit}
 	Alexander McCleary and Amit Patel.
 	\newblock Edit distance and persistence diagrams over lattices.
 	\newblock {\em arXiv preprint arXiv:2010.07337}, 2020.
 	
 	\bibitem[Mun17]{Munch:2017}
 	Elizabeth Munch.
 	\newblock A user's guide to topological data analysis.
 	\newblock {\em Journal of Learning Analytics}, 4(2):47--61, 2017.
 	
 	\bibitem[Pat18]{Patel:2018}
 	Amit Patel.
 	\newblock Generalized persistence diagrams.
 	\newblock {\em J. Appl. Comput. Topol.}, 1(3-4):397--419, 2018.
 	
 	\bibitem[Pes86]{pestov1986free}
 	V.~G. Pestov.
 	\newblock Free {B}anach spaces and representations of topological groups.
 	\newblock {\em Funktsional. Anal. i Prilozhen.}, 20(1):81--82, 1986.
 	
 	\bibitem[Rie16]{Riehl:2016}
 	Emily Riehl.
 	\newblock {\em Category theory in context}.
 	\newblock Aurora: Modern Math Originals. Dover Publications Inc., 2016.
 	
 	\bibitem[RR98]{rachev1998mass}
 	Svetlozar~T. Rachev and Ludger R\"{u}schendorf.
 	\newblock {\em Mass transportation problems. {V}ol. {I}}.
 	\newblock Probability and its Applications (New York). Springer-Verlag, New
 	York, 1998.
 	\newblock Theory.
 	
 	\bibitem[RT17]{Robinson:2017}
 	Andrew Robinson and Katharine Turner.
 	\newblock Hypothesis testing for topological data analysis.
 	\newblock {\em J. Appl. Comput. Topol.}, 1(2):241--261, 2017.
 	
 	\bibitem[SDB16]{Seversky_2016_CVPR_Workshops}
 	Lee~M. Seversky, Shelby Davis, and Matthew Berger.
 	\newblock On time-series topological data analysis: New data and opportunities.
 	\newblock In {\em Proceedings of the IEEE Conference on Computer Vision and
 		Pattern Recognition (CVPR) Workshops}, June 2016.
 	
 	\bibitem[ST20]{Skraba:2020}
 	Primoz Skraba and Katharine Turner.
 	\newblock Wasserstein stability for persistence diagrams.
 	\newblock 06 2020.
 	\newblock arXiv:2006.16824 [math.AT].
 	
 	\bibitem[Vil03]{villani2003topics}
 	C{\'e}dric Villani.
 	\newblock {\em Topics in optimal transportation}.
 	\newblock Number~58. American Mathematical Soc., 2003.
 	
 	\bibitem[Wea18]{weaver2018lipschitz}
 	Nik Weaver.
 	\newblock {\em Lipschitz algebras}.
 	\newblock World Scientific Publishing Co. Pte. Ltd., Hackensack, NJ, 2018.
 	\newblock Second edition of [ MR1832645].
 	
 	\bibitem[Wei13]{Weibel:kbook}
 	Charles~A. Weibel.
 	\newblock {\em The {$K$}-book}, volume 145 of {\em Graduate Studies in
 		Mathematics}.
 	\newblock American Mathematical Society, Providence, RI, 2013.
 	
 \end{thebibliography}

\end{document}